\documentclass{amsart}

\usepackage{pinlabel}  
\usepackage{graphicx}  
\usepackage{amscd}
\usepackage{amsmath}
\usepackage{bm}
\usepackage{enumerate}
\usepackage[dvipdfmx]{hyperref}
\usepackage{mathdots}
\usepackage{amssymb}

\title{Period matrices of some hyperelliptic Riemann surfaces}

\author{Yoshihiko Shinomiya}
\address{Mathematics Education, 
Faculty of Education College of Education, 
Academic Institute, Shizuoka University 
836 Ohya, Suruga-ku, Shizuoka 422-8529, JAPAN}
\email{shinomiya.yoshihiko@shizuoka.ac.jp}
\urladdr{}


\newtheorem{theorem}{Theorem}[section]    

\newtheorem{proposition}[theorem]{Proposition}   

\newtheorem{lemma}[theorem]{Lemma}          
\theoremstyle{definition}
\newtheorem{definition}[theorem]{Definition}    
\newtheorem*{remark}{Remark}             
\newtheorem{example}[theorem]{Example}   

\theoremstyle{plain}
 \newtheorem{corollary}[theorem]{Corollary}
 \newtheorem{notation}[theorem]{Notation}    

\newcommand{\thm}[1]{\begin{theorem} #1\end{theorem}}
\newcommand{\prop}[1]{\begin{proposition} #1\end{proposition}}
\newcommand{\lem}[1]{\begin{lemma} #1\end{lemma}}
\newcommand{\rem}[1]{\begin{remark} #1\end{remark}}
\newcommand{\pf}[1]{\begin{proof} #1\end{proof}}
\newcommand{\defi}[1]{\begin{definition} #1\end{definition}}
\newcommand{\exam}[1]{\begin{example} #1\end{example}}
\newcommand{\cor}[1]{\begin{corollary} #1\end{corollary}}
\newcommand{\note}[1]{\begin{notation} #1\end{notation}}
\newcommand{\dis}{\displaystyle}

\newcommand{\zero}{{\mathrm{Zero}}}

\newcommand{\interior}[1]{{\mathrm{Int}}(#1)}

\renewcommand{\zero}[1]{\mathrm{Zero}(#1)}
\newcommand{\jac}{\mathrm{Jac}}

\begin{document}

\begin{abstract}    
In this paper, we calculate period matrices 
of algebraic curves defined by 
$w^2=z(z^2-1)(z^2-a_1^2)(z^2-a_2^2)\cdots (z^2-a_{g-1}^2)$
for any $g\geq 2$ and  $a_1, a_2, \dots, a_{g-1}\in \mathbb{R}$ 
with $1<a_1<a_2<\cdots <a_{g-1}$. 
We construct these algebraic curves from Euclidean polygons.
A symplectic basis of these curves are 
given from the polygons.  
\end{abstract}

\maketitle

%
\section{Introduction}

Let $X$ be 
a compact Riemann surface 
or a smooth algebraic curve 
over $\mathbb{C}$. 
Assume that the surface $X$ is of genus $g \geq 2$.

\defi{[Symplectic Basis]
A basis 
$\left\{ \alpha_1, \dots, \alpha_g,  \beta_1, \dots, \beta_g \right\}$  
of $H_1(X, \mathbb{Z})$
 is called symplectic if  
$i(\alpha_k, \alpha_j)= i(\beta_k, \beta_j)=0$ and
$i(\alpha_k, \beta_j)=\delta_{k,j}$
hold 
for all $j$, $k$.
Here, $i(\cdot, \cdot)$ is the intersection number function.
}

Let $\Omega(X)$ be the space of holomorphic $1$-forms on $X$. 
Given a symplectic basis 
$\left\{ \alpha_1, \dots, \alpha_g,  \beta_1, \dots, \beta_g \right\}$  
of $H_1(X, \mathbb{Z})$ and 
a basis $\left\{\omega_1, \omega_2, \dots, \omega_g\right\}$ of $\Omega(X)$.
Define two square matrices of size $g$ by 
$A=\left[\int_{\alpha_j} \omega_k  \right]$
and 
$B=\left[\int_{\beta_j} \omega_k  \right]$.

\defi{[Period Matrix]
The matrix $\Pi=A^{-1}B$ is called the 
period matrix of $X$
 for $\left\{ \alpha_1, \dots, \alpha_g,  \beta_1, \dots, \beta_g \right\}$.
}
The period matrix depends only on the choice of 
symplectic basis of $H_1(X, \mathbb{Z})$.  
It is known that $\Pi$ is symmetric 
and its imaginary part $\mathrm{Im}(\Pi)$ 
is positive definite. 
That is, 
$\Pi$ is an element of the 
Siegel upper half-space  
$\mathcal{H}_g$ 
of degree $g$.  
Let $\left\{ \alpha_1^\prime, \dots, \alpha_g^\prime,  
\beta_1^\prime, \dots, \beta_g^\prime \right\}$  
be another symplectic basis of $H_1(X, \mathbb{Z})$. 
Then, there exists 
$
 T=\left[
\begin{array}{cc} 
P & Q\\
R & S
\end{array}
\right]
\in\mathrm{Sp}_{2g} (\mathbb{Z}) 
$
such that 
$
[\alpha_1^\prime,\dots, \beta_g^\prime ]
=[\alpha_1, \dots, \beta_g]T
$ 
holds. 
Here, $P$, $Q$, $R$, $S$ 
are square matrices of size $g$.  
Then, the period matrix $\Pi^\prime$ of $X$ 
for  $\left\{ \alpha_1^\prime, \dots, \alpha_g^\prime,  
\beta_1^\prime, \dots, \beta_g^\prime \right\}$  
is described as 
\begin{align*}
\Pi^\prime= (P\Pi +Q)(R \Pi +S)^{-1} .
\end{align*} 
This induces a map $\Phi$ from 
the moduli space  $\mathcal{M}_g$
of compact Riemann surfaces of genus $g$ 
to $\mathrm{Sp}_{2g} (\mathbb{Z}) \setminus \mathcal{H}_g$ 
which is called 
the Siegel modular variety of degree $g$.
The  Siegel modular variety 
is 
 the moduli space of
$g$-dimensional principally polarized abelian varieties.
For each equivalence class $[X]\in \mathcal{M}_g$, 
$\Phi([X])$ is the Jacobian variety $J(X)$ 
of $X$ which is 
defined 
by 
$\jac(X)=\mathbb{C}^g/\left(\mathbb{Z}^g+ \Pi \mathbb{Z}^g\right)$. 
The Torelli theorem states that 
two Riemann surfaces $X$ and $Y$ 
is conformal equivalent 
if and only if 
their Jacobian varieties $J(X)$ and $J(Y)$ 
are isomorphic as
polarized abelian varieties.
This implies that the map $\Phi$ is injective.
Therefore, period matrices are important data 
to study complex structures of Riemann surfaces. 
However, there are few examples of period matrices. 
It is difficult to find symplectic bases of Riemann surfaces in general.


For low genus case, 
the period matrix of 
the algebraic curve defined by 
$w^7=z(1-z)$ 
is calculated 
in \cite{TreTre84} and \cite{Tadokoro08}. 
This curve is of genus $3$. 
The Klein quartic curve 
is a curve in $\mathbb{C}P^3$ 
defined by $X Y^3+Y Z^3+Z X^3=0$. 
This curve is of genus $3$. 
The period matrix of this  curve  is 
 calculated in 
\cite{BraNor10},
\cite{Kamata02-1},
\cite{RauLew70},
\cite{RodRubGonVic97},
\cite{Schindler91},
\cite{Tadokoro08},
\cite{Yoshida99}
and 
\cite{Tadokoro08}.
Berry and Tretkoff \cite{BerTre92} 
calculated explicitly the period  matrix of 
Macbeath's curve 
which is of genus $7$. 
Kuusalo and N\"{a}\"{a}t\"{a}nen \cite{KuuNaa95}
calculated explicitly the period  matrices of 
algebraic curves defined by 
$w^2=z(z^4-1)$, 
$w^2=z^5-1$ 
and 
$w^2=z^6-1$. 

For generic genus, 
Schindler\cite{Schindler93} calculated 
the period matrices of the algebraic curves 
defined by 
$w^2=z^{2g+2}-1$, 
$w^2=z(z^{2g+1}-1)$ and  
$w^2=z(z^{2g}-1)$ 
for $g \geq 2$. 
Tashiro, Yamazaki, Ito 
and Higuchi \cite{TasYamItoHig96}
calculated the  period matrix 
of the algebraic curve 
defined by
$w^2=z^{2g+1}-1$. 
The explicit form of this period matrix 
is given by Tadokoro \cite{Tadokoro08}. 
Bujalance, Costa, Gamboa and Riera \cite{BujCosGamRie00}  
calculated period matrices of algebraic curves 
defined by 
$w^2=z^{2g+2}-1$ 
and 
$w^{2g+2}=z(z-1)^{g-1}(z+1)^{g+2}$. 
We know only these examples of period matrices 
for generic genus. 

In this paper, we calculate period matrices 
of algebraic curves defined by 
$$w^2=z(z^2-1)(z^2-a_1^2)(z^2-a_2^2)\cdots (z^2-a_{g-1}^2)$$ 
for any $g\geq 2$ and  $a_1, a_2, \dots, a_{g-1}\in \mathbb{R}$ 
with $1<a_1<a_2<\cdots <a_{g-1}$. 
In section \ref{section_construction}, 
we construct  Riemann surfaces from Euclidean polygons 
and show that these Riemann surfaces are the algebraic curves 
defined by the equations as above. 
This is done by finding their automorphisms. 
We also show that 
all algebraic curves 
defined by the equations as above 
are obtained by our construction from  Euclidean polygons. 
In section \ref{section_calculate}, 
we give symplectic bases of our algebraic curves 
and calculate the period matrices of them. 
In section \ref{section_example}, 
we give some examples of calculations of period matrices. 
Especially, for genus two case, we give period matrices explicitly. 
In section \ref{section_appendix}, we show that our algebraic curves are 
different from Schindler's four curves.

\section*{acknowledgement}
This work was supported by JSPS KAKENHI Grant Number 17K14184.
I would like to thank Masanori Amano 
for his careful reading and important comments. 
I am deeply grateful to Yuuki Tadokoro 
for leading the author to this area and
giving many valuable comments.

\label{section_introduction}
\section{Construction of Riemann Surfaces and their algebraic equations} \label{section_construction}

In this section, 
we construct hyperelliptic Riemann surfaces 
which we calculate their period matrices.
The Riemann surfaces are constructed from  
some rectangles. 
We describe them as algebraic curves over $\mathbb{C}$.

\subsection{Construction of Riemann Surfaces}\label{construction_of_R_s}
\label{construction_subsection}
Let $g \geq 2$ be a natural number.
We differentiate the construction by parity.
For a rectangle $P \subset \mathbb{C}$ 
with a horizontal side, 
we denote by $h(P)$ and $w(P)$ the height and
width of $P$, respectively.

\begin{figure}[h]
\labellist
\hair 0pt
\pinlabel $P_0$  at  45 45   
\pinlabel $P_1$  at  167 45 
\pinlabel $P_2$  at   167 120
\pinlabel $P_3$  at    288 120 
\pinlabel $P_{g-2}$  at   389 286
\pinlabel $P_{g-1}$  at    502 286
\endlabellist
\centering
\includegraphics[scale=0.35]{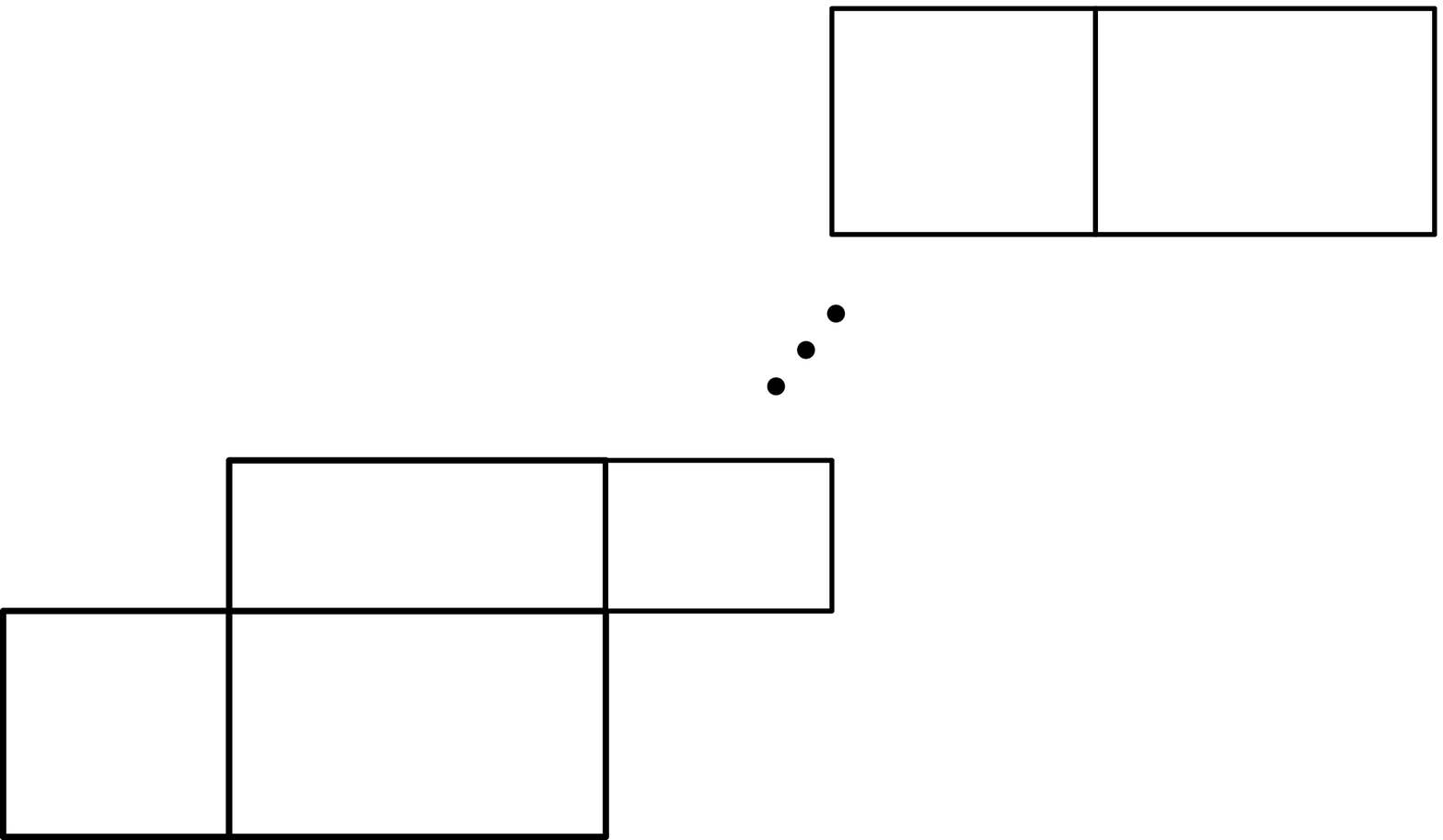} 
\caption{Rectangles $P_0, \dots, P_{g-1}$ if $g$ is even.}
\label{construction_even}
\end{figure}

Assume that $g$ is even.  
Let $P_0, \dots, P_{g-1}$ 
be rectangles satisfying the following:
\begin{itemize}
\item $P_i$ has a horizontal side for each $i=0, \dots, g-1$,
\item $P_0$ is a square,
\item $h(P_{2i})=h(P_{2i+1})$ for all $i=0, 1, \dots, \frac{g-2}{2}$, 
\item $P_{2i}$ and $P_{2i+1}$ have common vertical side and 
$P_{2i}$ is on the left of $P_{2i+1}$  for all $i=0, 1, \dots, \frac{g-2}{2}$, 
\item if $g\geq 4$, 
$w(P_{2i-1})=w(P_{2i})$ for all $i=1, 2, \dots, \frac{g-2}{2}$,
\item if $g\geq 4$, 
$P_{2i-1}$ and $P_{2i}$ have a common horizontal side and
$P_{2i-1}$ is below $P_{2i}$ 
for all $i=1, 2, \dots, \frac{g-2}{2}$
\end{itemize}
(see Figure \ref{construction_even}).

\begin{figure}[h]
\labellist
\hair 0pt
\pinlabel $P_0$  at  45 45   
\pinlabel $P_1$  at  45 114
\pinlabel $P_2$  at  154 114
\pinlabel $P_3$  at  154 161 
\pinlabel $P_{g-2}$  at  333 224
\pinlabel $P_{g-1}$  at  333 276
\endlabellist
\centering
\includegraphics[scale=0.35]{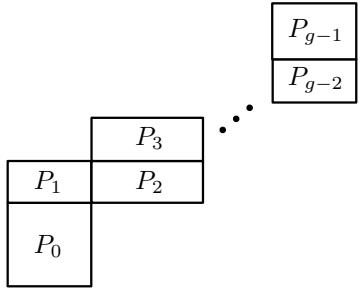} 
\caption{Rectangles $P_0, \dots, P_{g-1}$ if $g$ is odd.}
\label{paper_construction_odd}
\end{figure}

Next, assume that $g$ is odd.  
Let $P_0, \dots, P_{g-1}$ 
be rectangles satisfying the following:
\begin{itemize}
\item $P_i$ has a horizontal side for all $i=0, \dots, g-1$,
\item $P_0$ is a square,
\item $w(P_{2i})=w(P_{2i+1})$ for all $i=0, 1, \dots, \frac{g-3}{2}$, 
\item $P_{2i}$ and $P_{2i+1}$ have a common horizontal side and 
$P_{2i}$ is below $P_{2i+1}$  for all $i=0, 1, \dots, \frac{g-3}{2}$, 
\item 
$h(P_{2i-1})=h(P_{2i})$ for all $i=1, 2, \dots, \frac{g-2}{2}$,
\item 
$P_{2i-1}$ and $P_{2i}$ have common vertical side and
$P_{2i-1}$ is on the left of $P_{2i}$ 
for all $i=1, 2, \dots, \frac{g-2}{2}$
\end{itemize}
(see Figure\ref{paper_construction_odd}).

Let $g \geq 2$ and $P_0, P_1, \dots, P_{g-1}$ as above.
Let $l$ be the line passing through 
the upper left vertex and
the lower right vertex  of $P_0$.
Denote by $r_l$ the reflection about the line $l$.
We set 
$Q_i=r_l(P_i)$ for $i=1,2,\dots, g-1$ and 
\begin{align*}
P=P_0\cup \bigcup_{i=1}^{g-1} \left(P_i\cup Q_i\right).
\end{align*}
We make pairs of sides of $P$.
Each side of $P$ is paired with its ``opposite'' side.
For example, if  $g \geq 4$ is even, 
the left side of $P_{2i}$ and the right side of $P_{2i+1}$ 
are pair for all $i=0, 1, \dots, \frac{g-2}{2}$
and the lower side of $P_{2i-1}$ and 
the upper side of $P_{2i}$ are 
pair for all $i=1, 2, \dots, \frac{g-2}{2}$.
Identifying all pairs of sides of $P$ by parallel translation, 
we obtain a closed surface $X$ of genus $g$.
We induce a complex structure of $X$ 
so that the polygon $P$ gives one of the chars.
Then, $X$ is a Riemann surface of genus $g$.

\subsection{Algebraic equations of the Riemann Surfaces $X$}
Let $g \geq 2$ and $P$, $X$ as in subsection \ref{construction_of_R_s}.
In this subsection, we give 
an algebraic equation of the Riemann surface $X$.
We use the following notations. 

\begin{figure}[h]
\labellist
\hair 0pt
\pinlabel $p_0$  at  93 317
\pinlabel $p_1$  at   203 317
\pinlabel $p_2$  at    203 379
\pinlabel $p_3$  at    330 379
\pinlabel $p_4$  at    330 420
\pinlabel $o$  at    415 415
\pinlabel $q_1$  at  130 200 
\pinlabel $q_2$  at  50 200
\pinlabel $q_3$  at  50  80
\pinlabel $q_4$  at    -16 80
\pinlabel $o^\prime$  at   -15 5 
\endlabellist
\centering
\includegraphics[scale=0.35]{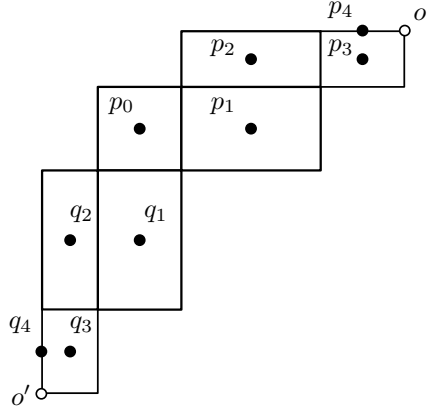} 
\caption{The points $p_0, \dots,p_{g}, q_1, \dots, q_{g}, o$ and $o^\prime$ if $g=4$.}
\label{points}
\end{figure}

\note{\label{notation} 
We name some points of $P$ and $X$ as follows 
(see Figure \ref{points});
\begin{itemize}
\item $o$ is the upper right vertex of $P_{g-1}$,
\item $o^\prime$ is the lower left vertex of $Q_{g-1}$,
\item $p_i$ is  the center of $P_j$ for $j=0, 1,\dots, g-1$,
\item $q_i$ is the center of $Q_j$ for $j=1,2,\dots, g-1$,
\item if $g$ is even, then $p_{g}$ is the midpoint of the upper side of $P_{g-1}$,
\item if $g$ is even, then $q_{g}$ is the midpoint of the left side of $Q_{g-1}$,
\item if $g$ is odd, then $p_{g}$ is the midpoint of the right side of $P_{g-1}$,
\item if $g$ is odd, then $q_{g}$ is the midpoint of the lower side of $Q_{g-1}$,
and
\item we use the same symbols to the points on $X$ corresponding to the above points.
\end{itemize}
Note that $o$ and $o^\prime$
coincide on $X$. 
We denote the point of $X$ by $o$.
}

Let $l$ be the line as in  subsection \ref{construction_of_R_s} 
and $r_l$ the reflection about $l$.
This reflection $r_l$ induces an antiautomorphism $\mu$ of $X$.
Moreover, we can construct an automorphism $\tau$ of $X$ 
from 
the $\frac{\pi}{2}$-rotation about $p_0$. 
Let $P^\prime$ be the image of $P$ by this rotation.
Then $P^\prime$ does not coincide with $P$.
However, we can reconstruct $P$ from $P^\prime$.
We cut $P^\prime$  along all sides 
which are shared by two of $P_0, P_1, \dots, P_{g-1}, Q_1, \dots, Q_{g-1}$
and glue the rectangles by the identification of sides of $P$.
Then the resulting polygon is $P$ with the identification of the opposite sides. 
Thus,  we obtain an automorphism $\tau$ of $X$.
For the maps $\tau$ and $\mu$, we have the following.

\prop{\label{automorphism}
The antiautomorphism $\mu$ is of order $2$. 
The automorphism $\tau$ is of order $4$. 
Moreover, $\mu^2$ is the hyperelliptic involution of $X$.
}

\pf{
The orders of $\tau$ and $\mu$ is clearly $2$ and $4$, respectively.
We show that $\mu^2$ is  the hyperelliptic involution of $X$.
By construction, $\mu^2$ maps each $P_i$ and $Q_j$ to itself.
Since $\mu^2$ acts on $P_i$ as $\pi$-rotation for each $i=0, 1, \dots, g-1$, 
$\mu^2$ fixes the center $p_i$ of $P_i$.
Moreover, $\mu^2$ fixes $p_{g}$.
By the same argument, $\mu^2$ fixes $q_1, q_2, \dots,q_{g-1}$ and $q_g$. 
Finally, $\mu^2$ fixes $o$.
Therefore, $\mu^2$ has $2g+2$ fixed points.
This means that $\mu^2$ is the hyperelliptic involution of $X$.
}

We now describe $X$ as an algebraic curve.

\thm{\label{algebraic}
The Riemann surface $X$ is 
the algebraic curve defined by 
\begin{align*}
w^2=z(z^2-1)(z^2-a_1^2)\dots (z^2-a_{g-1}^2).
\end{align*}
Here, $a_1, \dots, a_{g-1}$ are real number satisfying 
$1<a_1< a_2<\cdots< a_{g-1}$. 
Moreover, we have 
$p_0=(0, 0), p_1=(1, 0), p_{j+1}=(a_j, 0), 
q_1=(-1, 0), q_{j+1}=(-a_j, 0)$ 
for all $j=1, \dots, g-1$. 
Let $(z, w)$ be a point of the algebraic curve $X$. 
Then, we have 
$\tau(z, w)=(-z, i w)$ and $\mu(z, w)=(-\bar{z}, -i\bar{w})$.
}

\rem{
Hereafter, we choose a branch of the square root 
so  that $\sqrt{-x}=-i\sqrt{x}$ holds for all $x>0$.
}

To prove Theorem \ref{algebraic}, 
we show the following lemma.

\lem{\label{integral}
Let $f(x)$ be a real valued function defined on an interval $I$. 
Assume that the function $\frac{1}{\sqrt{f(x)}}  : I \to \hat{\mathbb{R}}$ is integrable on $I$. 
Then, we have
\begin{align*}
\int_I \dfrac{d x}{\sqrt{-f(x)}} 
=i \cdot \overline{\int_I \frac{d x}{\sqrt{f(x)}}}.
\end{align*}
}

\pf{
Set $I_1=\left\{x\in I : f(x) \geq 0 \right\}$ 
and 
$I_2=\left\{x\in I : f(x) < 0 \right\}$.
Then the integral 
$\dis \int_{I_1} \frac{d x}{\sqrt{f(x)}}$ 
is a real number 
and
$\dis \int_{I_2} \frac{d x}{\sqrt{f(x)}}$
is a  pure imaginary number. 
Moreover, the equation 
\begin{align*}
-i \int_{I_2} \frac{d x}{\sqrt{f(x)}}
= -i \int_{I_2} \frac{d x}{-i \sqrt{-f(x)}}
=\int_{I_2} \frac{d x}{ \sqrt{-f(x)}}
\end{align*}
holds. 
Therefore, we have 
\begin{align*}
\int_I \dfrac{dx}{\sqrt{-f(x)}} 
&=\int_{I_1} \dfrac{dx}{\sqrt{-f(x)}} +\int_{I_2} \dfrac{dx}{\sqrt{-f(x)}} \\
&=\int_{I_1} \dfrac{dx}{-i \sqrt{f(x)}} -i \int_{I_2} \frac{dx}{\sqrt{f(x)}}\\
&=i \int_{I_1} \dfrac{dx}{ \sqrt{f(x)}} -i \int_{I_2} \frac{dx}{\sqrt{f(x)}}\\
&=i \left( \int_{I_1} \dfrac{dx}{ \sqrt{f(x)}} - \int_{I_2} \frac{dx}{\sqrt{f(x)}} \right)\\
&=i \cdot \overline{\int_I \frac{dx}{\sqrt{f(x)}}}.  \qedhere 
\end{align*}
}

\pf{[Proof of Theorem \ref{algebraic}]
By Proposition \ref{automorphism}, 
$\tau^2$ is the hyperelliptic involution of $X$
and $Y=X/\langle \tau^2 \rangle$ is the Riemann sphere. 
Let $\varphi : X \to Y=\hat{\mathbb{C}}$ be the natural projection. 
We may assume that $\varphi(p_0)=0$, $\varphi(p_1)=1$ and $\varphi(o)=\infty$. 
The automorphism $\tau$ of $X$ induces an automorphism $\tau^\prime$ of $Y$. 
The automorphism $\tau^\prime$ is of order $2$ and 
fixes $\varphi(p_0)=0$ and $\varphi(o)=\infty$. 
Hence, $\tau^\prime$ is the M\"{o}bius transformation which is of the form $\tau^\prime(z)=-z$. 
From this, we have $\varphi(p_j)=-\varphi(q_j)$ for all $j=1, 2, \dots, g$. 
Since $\mu \tau^2=\tau^2 \mu$ holds, 
the antiautomorphism $\mu$ of $X$ induces an automorphism $\mu^\prime$ of $Y$. 
The antiautomorphism $\mu^\prime$ is of order $2$ and 
fixes $\varphi(p_0)=0$ and $\varphi(o)=\infty$. 
Moreover, $\mu^\prime$ maps $\varphi(p_1)=1$ to $\varphi(q_1)=-1$. 
Therefore,  $\tau^\prime$ is the antiautomorphism 
which is of the form $\mu^\prime(z)=-\bar{z}$. 
For each $j=2, 3, \dots, g-1$, we have
\begin{align*}
\varphi(p_j)
=-\varphi(q_j) 
=-\varphi(\mu(p_j))
=-\mu^\prime(\varphi(p_j))
=\overline{ \varphi(p_j)}.
\end{align*}
This implies that $\varphi(p_j)$ is a real number 
for each  $j=2, 3, \dots, g-1$.
We set $a_{j-1}=\varphi(p_j)$ 
for each $j=2,3, \dots, g-1$. 
Now,  $X$ is described by the algebraic equation
\begin{align*}
 w^2=z(z^2-1)(z^2-a_1^2)(z^2-a_2^2) \cdots(z^2-a_{g-1}^2).
\end{align*}

Next, we show that $1<a_1<a_2<\cdots <a_{g-1}$ holds. 
Set $f(z)=z(z^2-1)(z^2-a_1^2)(z^2-a_2^2) \cdots(z^2-a_{g-1}^2)$.
Assume that the polygon $P$ is in the $z$-plane.
Let $\omega$ be the holomorphic $1$-form on $X$ 
induced by the $1$-form $dz$ on the polygon $P$.
The $1$-form $\omega$ has a unique zero $o$ of order $2g-2$. 
The holomorphic $1$-form $\frac{dz}{w}$ is also 
has a unique zero $o$ since $o$ corresponds to $\infty$ via $\varphi$. 
There exists a constant $c \in \mathbb{C}^\ast$ such that $\omega=c\frac{dz}{w}$ holds.
Let $\zeta$ be a local coordinate system of $X$ induced by $\interior{P}$. 
We may assume that $\zeta(p_0)=0$. 
We describe points $p\in X$  by $(z,w) \in \mathbb{C}^2$ 
satisfying the above equation. 
Then, for sufficiently small neighborhood $U$ of $p_0$,  
the local coordinate system $\zeta$ is represented by
\begin{align*}
\zeta(p)
=\int_{p_0}^p \omega
=\int_{(0,0)}^{(z,w)} c\frac{dz}{w}
=c\int_{0}^{z} \frac{dz}{\sqrt{f(z)}}.
\end{align*}
Let $C$ be a lift of the real axis via the projection $\varphi : X \to X/\langle \tau^2\rangle$. 
Let $F: C \to \mathbb{C}$  be the analytic continuation of $\zeta$ along $C$. 
If $g$ is even,  we have 
\begin{align*}
|p_0 p_1|
=\int_{p_0}^{p_1}  \omega
= c\int_{0}^{1} \frac{dz}{\sqrt{f(z)}}.
\end{align*}
Moreover, by Lemma \ref{integral}, 
\begin{align*}
|p_0p_1|=
|p_0 q_1|
&=i\int_{p_0}^{q_1} \omega 
=ic\int_0^{-1} \frac{dz}{\sqrt{f(z)}}
=ic\int_0^{1} \frac{-dz}{\sqrt{-f(x)}}
=c \cdot \overline{\int_{0}^{1} \frac{dz}{\sqrt{f(z)}}. }
\end{align*}
Therefore, $\dis \int_{0}^{1} \frac{dz}{\sqrt{f(z)}}$ and $c$ are real numbers. 
If there exists $l\in \left\{1, 2, \dots, g-1\right\}$ 
such that $0<|a_l|<1$, 
the integral 
$\dis \int_{0}^{1} \frac{dz}{\sqrt{f(z)}}$
is not real. 
Thus, $1<|a_l|$ holds for all $l\in \left\{1, 2, \dots, g-1\right\}$ . 
Let $C_0$ be the path from $p_0$ to $p_1$ along $C$.
Then $F(C)$ is a horizontal segment $p_0p_1$ in $P$.
Moreover, the image $F(C)$  consists of horizontal segments and vertical segments. 
This implies that $F(C)$ pass through $p_0, p_1, \dots, p_g, o$ in this order. 
Therefore, we have $1<a_1<a_2<\cdots <a_{g-1}$. 
If $g$ is odd, by the same way, we can see that 
$\dis \int_{0}^{1} \frac{dz}{\sqrt{f(z)}}$ is a pure imaginary number, $c$ is a real number 
and hence, $1<a_1<a_2<\cdots <a_{g-1}$ holds. 

Finally, we show $\tau(z, w)=(-z, iw)$ and $\mu(z, w)=(-\bar{z}, -i\bar{w})$ for all $(z, w)\in X$. 
Since $\tau^\prime(z)=-z$ holds, the first coordinate of $\tau(z, w)$ is $-z$. 
The equation $f(-z)=-f(z)=-w^2$ implies that 
the second coordinate of $\tau(z, w)$ is $\pm i w$. 
We set  $\tau(z, w)=(-z, \varepsilon i w)$.  
Then, 
$
\int_{\tau(C_0)} \omega
 =i |p_0 q_1|
 =i |p_0 p_1|
$ 
and 
\begin{align*}
\int_{\tau(C_0)} \omega
= \int_{\tau(C_0)} c\frac{dz}{w}
=\int_{C_0} \tau^\ast \left( c\frac{dz}{w} \right)
=\int_{C_0}  c\frac{-dz}{i \varepsilon w} 
=\frac{i}{\varepsilon} \int_{C_0}  c\frac{dz}{ w} 
=\frac{i}{\varepsilon} |p_0 p_1|
\end{align*}
hold. 
Hence, we have $\varepsilon =1$ 
and  $\tau(z, w)=(-z, iw)$ for all $(z, w)\in X$. 
}

\subsection{Representation of Algebraic Curves by Polygon}

In subsection \ref{construction_of_R_s}, 
we show that the Riemann surface
which is obtained from the polygon $P$ 
is the algebraic curve defined by 
\begin{align*}
 w^2=z(z^2-1)(z^2-a_1^2)(z^2-a_2^2) \cdots(z^2-a_{g-1}^2) 
\end{align*}
for some $a_1, a_2, \dots, a_{g-1}$ ($1<a_1< a_2<\cdots< a_{g-1}$).
In this subsection, we show the following. 

\thm{\label{representation}
Let $g \geq 2$.  
For all  $a_1, a_2, \dots, a_{g-1}$ ($1<a_1< a_2<\cdots< a_{g-1}$), 
the algebraic curve 
\begin{align*}
 w^2=z(z^2-1)(z^2-a_1^2)(z^2-a_2^2) \cdots(z^2-a_{g-1}^2) 
\end{align*} 
is obtained from the polygon 
\begin{align*}
P=P_0\cup \bigcup_{i=1}^{g-1} \left(P_i\cup Q_i\right)
\end{align*} 
as in 
Figure \ref{construction_even} or Figure \ref{paper_construction_odd} 
by adjusting the lengths of rectangles. 
Here, $P_0$ is a square.
}

We prove this theorem by the theory of translation surfaces.
Hereafter, we assume that $X$ is 
a compact Riemann surface of genus $g \geq 2$.

\defi{[Translation Surface] \label{translation_surface}
A  translation surface $(X, \omega)$ 
is a pair of a Riemann  surface $X$ and 
a non-zero holomorphic $1$-form $\omega$ on $X$.
The zeros of $\omega$ are called singular points  
of the  translation surface $(X, \omega)$. 
Denote by $\zero{\omega} $
the set of all zeros of $\omega$.
}

Let 
$\omega$ be  a non-zero holomorphic $1$-form $\omega$ on $X$.
If $p_0 \in X$ is not a zero of $\omega$, there is a  neighborhood $U$ 
such that 
\begin{eqnarray*}
z=\int_{p_0}^p \omega : U \rightarrow \mathbb{C}
\end{eqnarray*}
is a chart of $X$. 
Let $(U, z)$ and $(V, w)$ be such charts with $U \cap V \not = \emptyset$.
The transition function is of the form $w= z+(\mbox{const.})$. 
Hence, $(X, \omega)$ is a surface with a Euclidean structure on $X-\zero{\omega}$.
If $p_0$ is a zero of $\omega$ of order $n$, 
there exists a chart $\zeta  : U \to \mathbb{C}$ of $X$ 
such that $\zeta(p_0)=0$ $\omega=\zeta^n d\zeta$. 
Then, 
$$z=\int_{p_0}^p \omega
=\frac{1}{n+1} \zeta^{n+1}
 : U \rightarrow \mathbb{C}$$
is a chart around $p_0$. 
With respect to this chart, 
the angle around $p_0$ is $2\pi (n+1)$.

We define some terminologies for translation surfaces. 
Let $(X, \omega)$ be a translation surface.
We consider geodesics
with respect to the Euclidean structure  on $(X, \omega)$.

\defi{[Saddle connection]
A saddle connection of  a translation surface $(X, \omega)$
is a geodesic segment on $(X, \omega)$ 
whose end points are singular points 
and containing no singular points in its interior. 
Note that the end points of a saddle connection 
may be same. 
}

For closed geodesics on translation surfaces, 
we have the following.

\prop{
Let $\gamma$ be a closed geodesic on a translation surface $(X, \omega)$.
Then, one of the following two holds; 
\begin{enumerate}
\item $\gamma$ contains no singular points, 
\item $\gamma$ is a  concatenation of saddle connections. 
Assume that  $\gamma$ passes  a saddle connection $s_2$ 
just after a saddle connection $s_1$ and the switch 
occurs at a singular point $p_0$. 
The angles at $p_0$ in both sides of $s_1 \cup s_2$ 
are greater than or equal to $\pi$. 
\end{enumerate}
}

\defi{
Let $\gamma$ be a closed geodesic on 
a translation surface $(X, \omega)$ 
which does not contain singular points. 
The geodesic $\gamma$ is horizontal (resp. vertical)
if it is horizontal (resp. vertical)
 with respect to the Euclidean structure of  $(X, \omega)$. 
}

We construct polygon $P$ as in Theorem \ref{representation} 
from horizontal and vertical closed geodesics. 
Then, we use maximal cylinders for the geodesics.

\defi{[Maximal Cylinder]\label{maximal}
Let $\gamma$ be a closed geodesic on 
a translation surface $(X, \omega)$ 
which does not contain singular points. 
The maximal cylinder $R_\gamma$ for $\gamma$ 
is the union of all 
closed geodesics on $(X, \omega)$ 
which do not contain singular points 
and are homotopic to $\gamma$.
}

By Definition \ref{maximal}, we have the following. 

\prop{
Let $\gamma$ be a closed geodesic on 
a translation surface $(X, \omega)$ 
which does not contain singular points. 
The maximal cylinder  $R_\gamma$ for $\gamma$ 
is a cylinder each of whose boundary components 
are closed geodesics constructed from 
saddle connections that are parallel to $\gamma$.
}

We now prove Theorem \ref{representation}.
Choose $a_1, a_2, \dots, a_{g-1}$ 
so that $1<a_1< a_2<\cdots< a_{g-1}$ 
and set 
$f(z)=
 w^2=z(z^2-1)(z^2-a_1^2)(z^2-a_2^2) \cdots(z^2-a_{g-1}^2)
$.
Let $X$ be the algebraic curve $w^2=f(z)$ 
and $\varphi : X \to \hat{\mathbb{C}} ; (z, w) \mapsto z$ 
the natural projection. 
Set $o=\varphi^{-1}(\infty)$. 
The algebraic curve $X$ has 
an automorphism $\tau : X\to X ; (z,w) \to (-z, iw)$
.
%
The map $\tau$ is of order $4$ and 
$\tau^2$ is the hyperelliptic involution of $X$. 
The fixed points of $\tau^2$ are 
$(0,0), (\pm 1,0), (\pm a_1, 0), (\pm a_2, 0), \dots, 
(\pm a_{g-1}, 0)$ and $o$.  
Let $\delta_0, \delta_1, \dots, \delta_{g-1}$ 
be preimages of the intervals 
$[0, 1], [1,a_1], [a_1, a_2], \dots, [a_{g-2}, a_{g-1}]$, respectively. 
Set $\delta_j^\ast=\tau(\delta_j)$ for each $j=0, 1, \dots, g-1$. 
Then $\delta_j$ and $\delta_j^\ast$ are simple closed curves on $X$.

Let $i(\cdot,\cdot)$ be the geometric intersection number function. 
By construction, we have the following proposition.

\prop{\label{intersection}
The following holds. 
\begin{enumerate}[(1)]
\item 
For any $j=0,1, \dots, g-2$, 
$i(\delta_j, \delta_{j+1})=1$
holds.  
The curves $\delta_j$ and $\delta_{j+1}$ 
intersect only at $(a_j, 0)$.
If 
$|j-k| \not =1$
then 
$i(\delta_j, \delta_k)=0$ 
holds.

\item 
For any $j=0,1, \dots, g-2$, 
$i(\delta_j^\ast, \delta_{j+1}^\ast)=1$
holds.  
The curves $\delta_j^\ast$ and $\delta_{j+1}^\ast$
intersect only at $(-a_j, 0)$.
If 
$|j-k| \not =1$
then 
$i(\delta_j^\ast, \delta_k^\ast)=0$
holds.

\item 
The equation $i(\delta_0,  \delta_0^\ast)=1$ holds. 
The curves $\delta_0$ and $\delta_0^\ast$
intersect only at $(0, 0)$.
If $(j, k) \not= (0,0)$, then 
$i(\delta_j, \delta_k^\ast)=0$ holds.
\end{enumerate}
}

Set $\omega_0=\dfrac{dz}{w}$. 
The holomorphic differential $\omega_0$ 
has a unique zero at $o$. 
Therefore, $o$ is the unique singular point of 
the translation surface $(X, \omega_0)$. 
The angle around $o$ is $2\pi(2g-1)$ on $(X, \omega_0)$.
Moreover, we have the following proposition. 

\prop{\label{cross}
The simple closed curves 
$\delta_0, \dots, \delta_{g-1}, \delta_0^\ast,\dots,  \delta_{g-1}^\ast$
are closed geodesics on the translation surface  $(X, \omega_0)$. 
If $g$ is even, then 
$\delta_0$, $\delta_2$, $\dots$, $\delta_{g-2}$, 
$\delta_1^\ast$, $\delta_3^\ast$, $\dots$, $\delta_{g-1}^\ast$ 
are horizontal and 
$\delta_1$, $\delta_3$, $\dots$, $\delta_{g-1}$,
$\delta_0^\ast$, $\delta_2^\ast$, $\dots$, $\delta_{g-2}^\ast$
are vertical on $(X, \omega_0)$. 
If $g$ is odd, then
$\delta_0$, $\delta_2$, $\dots$, $\delta_{g-2}$, 
$\delta_1^\ast$, $\delta_3^\ast$, $\dots$, $\delta_{g-1}^\ast$ 
are vertical and 
$\delta_1$, $\delta_3$, $\dots$, $\delta_{g-1}$, 
$\delta_0^\ast$, $\delta_2^\ast$, $\dots$, $\delta_{g-2}^\ast$
are horizontal  on $(X, \omega_0)$.
}

\pf{
We prove the case if $g$ is even. 
Since $f(t)> 0$ holds for any $t\in (0, 1)$, 
$\int_0^t \frac{dz}{\sqrt{f(z)}}>0$ holds. 
Hence, $\delta_0$ is a horizontal closed geodesic on  $(X, \omega_0)$.
Since $f(t)< 0$ holds for any $t\in (1, a_1)$, 
$\int_1^t \frac{dz}{\sqrt{f(z)}} \in i\mathbb{R}$ holds. 
Hence, $\delta_1$ is a vertical closed geodesic on  $(X, \omega_0)$.
We can prove the others by the same way.
}

Let $\nu: X \to X;  (z,w) \mapsto  (\bar{z}, (-1)^g \bar{w})$ 
be an antiautomorphism  of $X$. 
Set $R_j=R_{\delta_j}$ and $R_j^\ast=R_{\delta_j^\ast}$ (see Definition \ref{maximal}). 
We regard each boundary component of the cylinder $R_j$ (resp. $R_j^\ast$)
as a closed curve which is homotopic to $\delta_j$ (resp. $\delta_j^\ast$).
We denote the boundary components 
by $\partial_1 R_j$ and $\partial_2 R_j$ 
(resp. $\partial_1 R_j^\ast$ and $\partial_2 R_j^\ast$).

\lem{\label{cylinder_inv}
For any $j = 0, 1, \dots, g-1$, 
we have the following; 
\begin{enumerate}[(1)]
\item $\tau(R_j)= R_j^\ast$ and  $\tau(R_j^\ast)= R_j$ holds.
\item \label{hyp_inv}
 The simple closed geodesics $\delta_j$ and  $\delta_j^\ast$
are invariant under $\tau^2$, respectively. 
Two boundary components $\partial_1 R_j$ and $\partial_2 R_j$
(resp. $\partial_1 R_j^\ast$ and $\partial_2 R_j^\ast$)
are permuted by $\tau^2$.
%

\item\label{nu_inv}  
Simple closed geodesics $\delta_j$ and $\delta_j^\ast$ 
are invariant under $\nu$, respectively. 
The  geodesic 
$\delta_j$ is pointwise fixed by $\nu$ if and only if $j$ is even.
Then, $\nu$ is represented as the reflection about  $\delta_j$ 
in a sufficiently small neighborhood of  $\delta_j$.
The  geodesic $\delta_j^\ast$ 
is pointwise fixed by $\nu$ if and only if $j$ is odd. 
Then, $\nu$ is represented as the reflection about  $\delta_j^\ast$ 
in a sufficiently small neighborhood of  $\delta_j^\ast$.

\item \label{nu_inv_cyl}
The cylinders $R_j, R_j^\ast$ are invariant under 
$\nu$, respectively. 
If $j$ is even, 
then two boundary components 
$\partial_1 R_j$ and $\partial_2 R_j$  
are permuted by $\nu$ and 
each boundary component of $R_j^\ast$ 
are invariant under  $\nu$. 
If $j$ is odd, 
then two boundary components  $\partial_1 R_j^\ast$ and $\partial_2 R_j^\ast$ 
are permuted by $\nu$ and 
each boundary component of $R_j$ 
are invariant under  $\nu$. 
\end{enumerate}
}

\pf{
\begin{enumerate}[(1)]
\item  Since $\delta_j^\ast=\tau(\delta_j)$, we obtain the claim. 
\item Since 
$\tau^\ast \dfrac{dz}{w}=i \dfrac{dz}{y}$ 
holds,
the automorphism 
$\tau$ acts on $(X, \omega_0)$ as $\frac{\pi}{2}$-rotation.
The closed geodesic $\delta_j$ passes through two fixed points of $\tau^2$. 
Hence, $\delta_j$ is invariant under $\tau^2$ and 
two boundary components of $R_j$ 
are permuted by $\tau^2$.
By the same argument as above, 
$\delta_j^\ast$ is invariant under $\tau^2$ and 
two boundary components of $R_j^\ast$ 
are permuted by $\tau^2$.

\item \label{nu_inv}
We prove the case where $g$ is even. 
The case where $g$ is odd is proved by the same way.
Assume that $g$ is even. 
Since 
$\nu^\ast \left(\dfrac{dz}{w} \right)= \overline{\left(\dfrac{dz}{w}\right)}$ 
holds, 
$\nu$ preserves horizontal slopes and vertical slopes of segments on $(X, \omega_0)$. 
The geodesic $\delta_j$ passes through 
$(a_{j-1}, 0)$ and $(a_j, 0)$ 
which are fixed points of $\nu$. 
Thus, $\delta_j$ is invariant under $\nu$.
By the same argument,  
$\delta_j^\ast$ is invariant under $\nu$.
Moreover, 
the set of all  fixed points of $\nu$ is
\begin{align*}
 \left\{ (z,w)\in X: z,w\in \mathbb{R} \right\}
=\delta_0\cup \delta_2 \cup \dots \cup \delta_{g-2}
\cup \delta_1^\ast \cup \delta_3^\ast \cup \dots \cup \delta_{g-1}^\ast .
\end{align*}
Therefore, we obtain the claim.

\item It is clear by (\ref{nu_inv}). 

\end{enumerate}
}

\lem{\label{twice}
For any $k=1, 2$ and $j = 0, 1, \dots, g-1$, 
the boundary component $\partial_k R_j$
passes through the singular point $o$ of $(X, \omega_0)$ 
at least twice.  
}

\pf{

By construction,  the curve $\partial_k R_j$ 
passed through $o$ at least once. 
Let $s^\prime$ be a perpendicular segment from $o$ to $\delta_j$ 
in $R_j$. 
If $j$ is even, we set $s=s^\prime \cup \nu(s^\prime)$.
If $j$ is odd,  we set 
$s=\tau^{-1}\left(\tau(s^\prime)\cup \nu\tau(s^\prime) \right)
=s^\prime \cup \tau^{-1}\nu \tau(s^\prime)$. 
I both cases, $s$ is a segment whose end points are $o$ 
and which is orthogonal to $\delta_j$. 
We show that $s$ is not invariant under $\tau^2$. 
If $s$ is invariant under $\tau^2$, 
then the midpoint of $s$, say $p$, is a fixed point of $\tau^2$. 
Thus $p=(a_{j-1}, 0)$ or $(a_j, 0)$. 
Here, we set $a_{-1}=0$, $a_0=1$. 
If $p=(a_{j-1}, 0)$, then $s$ is contained in $\delta_{j-1}$. 
This contradicts that $\delta_{j-1}$ does not pass through $o$. 
If $p=(a_{j}, 0)$, then $s$ is contained in $\delta_{j+1}$. 
This contradicts that $\delta_{j+1}$ does not pass through $o$. 
Therefore, $s$ is not invariant under $\tau^2$ and 
$\tau^2(s) \neq s$ holds. 
Now, $\tau^2(s)$ is a segment connecting 
two boundaries $\partial_1 R_j$ and $\partial_2 R_j$  in $R_j$ 
and whose end points are $o$.
Thus, we obtain the claim.
}

\pf{[Proof of Theorem \ref{representation}]
Let $c_j$ be the number of times that 
the boundary component $\partial_1 R_j$ passes through 
the singular point $o$ of $(X, \omega_0)$ 
for each $j=0,1,\dots, g-1$.
Then 
$\partial_2 R_j, 
\partial_1 R_j^\ast$ and $\partial_2 R_j^\ast
$
pass through $o$ just $c_j$ times, respectively. 
Set $L_j=\interior{R_j} \cap \interior{R_{j+1}}$  for $j=0, 1, \dots, g-2$.
Then $L_j$ is a rectangle  in $(X, \omega_0)$.
By Lemma \ref{cylinder_inv} (\ref{hyp_inv})  and (\ref{nu_inv_cyl}), 
$L_j$ is invariant under $\tau^2$ and $\nu$. 
This implies that if one of the vertices of $L_j$ is $o$, then 
all vertices are $o$.  
It is also true for rectangles $L=\interior{R_0}\cap  \interior{R_0^\ast}$ and 
$L_j^\ast=\interior{R_j^\ast}\cap  \interior{R_{j+1}^\ast}$ for $j=0, 1, \dots, g-2$.
Assume that $d$ of $L, L_0, \dots, L_{g-2}, L_0^\ast, \dots, L_{g-2}^\ast$
have $o$ as their vertices.
Since the angle around $o$ in $(X, \omega_0)$ is $(4g-2)\pi$, 
the inequality 
\begin{align*}
 2\sum_{j=0}^{g-1}2\pi c_j -2\pi d \leq (4g-2)\pi
\end{align*}
holds.
Since $c_j \geq 2$ for all 
$j= 0, 1, \dots, g-2$
 by Lemma \ref{twice}, we have 
\begin{align*}
2g-3+2c_{g-1} \leq d.
\end{align*}
Since  $d\leq 2g-1$ by definition,   
we obtain  
$c_{g-1}=1$,  $d=2g-1$ and 
$c_j = 2$
for all $j=0, 1, \dots, g-2$.
This implies that 
$\partial_k R_j$ 
(resp. $\partial_k R_j^\ast$ ) 
passes through 
the singular point $o$ just twice
for all $j=1, 2, \dots, g-2$ 
and $k=1, 2$ 
and  all vertices of rectangles 
$L, L_0, \dots, L_{g-2}, L_0^\ast, \dots, L_{g-2}^\ast$
are the singular point $o$. 
Therefore, we have
\begin{align*}
 X=\overline{L}\cup \bigcup_{j=0}^{g-2} \overline{L_j} \cup \bigcup_{j=0}^{g-2} \overline{L_j^\ast }.
\end{align*}
Moreover, since $L$ is invariant under $\tau$, 
$L$ is a square. 
}

\cor{
Let $g \geq 2$.
Let $a_1, a_2, \dots, a_{g-1}\in \mathbb{R} $ 
with $1<a_1<a_2<\cdots<a_{g-1}$ 
and set 
$a_{-1}=0$, $a_0=1$, $a_g=\infty$
and 
$f(z)=z(z^2-1)(z^2-a_1^2)\cdots(z^2-a_{g-1}^2)$. 
We also set 
\begin{align*}
I_j=\int_{a_{j-1}}^{a_j} \dfrac{dz}{\sqrt{|f(z)|}}
\end{align*}
for $j=0,1, \dots, g$.
Then the equation 
\begin{align*}
\sum_{j=0}^{g} (-1)^{[\frac{j+1}{2}]} I_j=0
\end{align*}
holds.
}

\pf{
Here, we prove the case where $g$ is even. 
By Theorem \ref{representation}, 
the algebraic curve $X$ defined by $w^2=f(z)$ is  
obtained from the polygon 
\begin{align*}
P=P_0\cup \bigcup_{i=1}^{g-1} \left(P_i\cup Q_i\right) 
\end{align*} 
as in 
Figure \ref{construction_even} 
by adjusting the lengths of rectangles. 
Here, $P_0$ is a square.
The length of horizontal edge of $L$ is represented by 
\begin{align*}
2\left( I_0-I_2+\cdots +(-1)^{[\frac{g-1}{2}]} I_{g-2}+(-1)^{[\frac{g+1}{2}]} I_g \right). 
\end{align*}
The length of vertical edge of $L$ is represented by 
\begin{align*}
2\left( I_1-I_3 +\dots+(-1)^{[\frac{g-4}{2}]} I_{g-3} +(-1)^{[\frac{g-2}{2}]} I_{g-1} \right)
\end{align*}
Since $L$ is a square, 
we have 
\begin{align*}
&I_0-I_2+\cdots +(-1)^{[\frac{g-1}{2}]} I_{g-2}+(-1)^{[\frac{g+1}{2}]} I_g\\
=&I_1-I_3 +\dots+(-1)^{[\frac{g-4}{2}]} I_{g-3} +(-1)^{[\frac{g-2}{2}]} I_{g-1}
\end{align*}
From this, we obtain the claim.
}

\section{Calculation of Period Matrices}\label{section_calculate}

\subsection{Construction of symplectic bases}\label{basis}
\begin{figure}[h]
\labellist
\hair 0pt
\pinlabel $\alpha_1$  at  -25 53
\pinlabel $\alpha_2$  at -25 144
\pinlabel $\alpha_3$  at 75 234
\pinlabel $\alpha_4$  at 160 323
\pinlabel $\beta_1$  at 585 210
\pinlabel $\beta_2$  at 675 305
\pinlabel $\beta_3$  at 765 390
\pinlabel $\beta_4$  at 855 390 
\endlabellist
\centering
\includegraphics[scale=0.3]{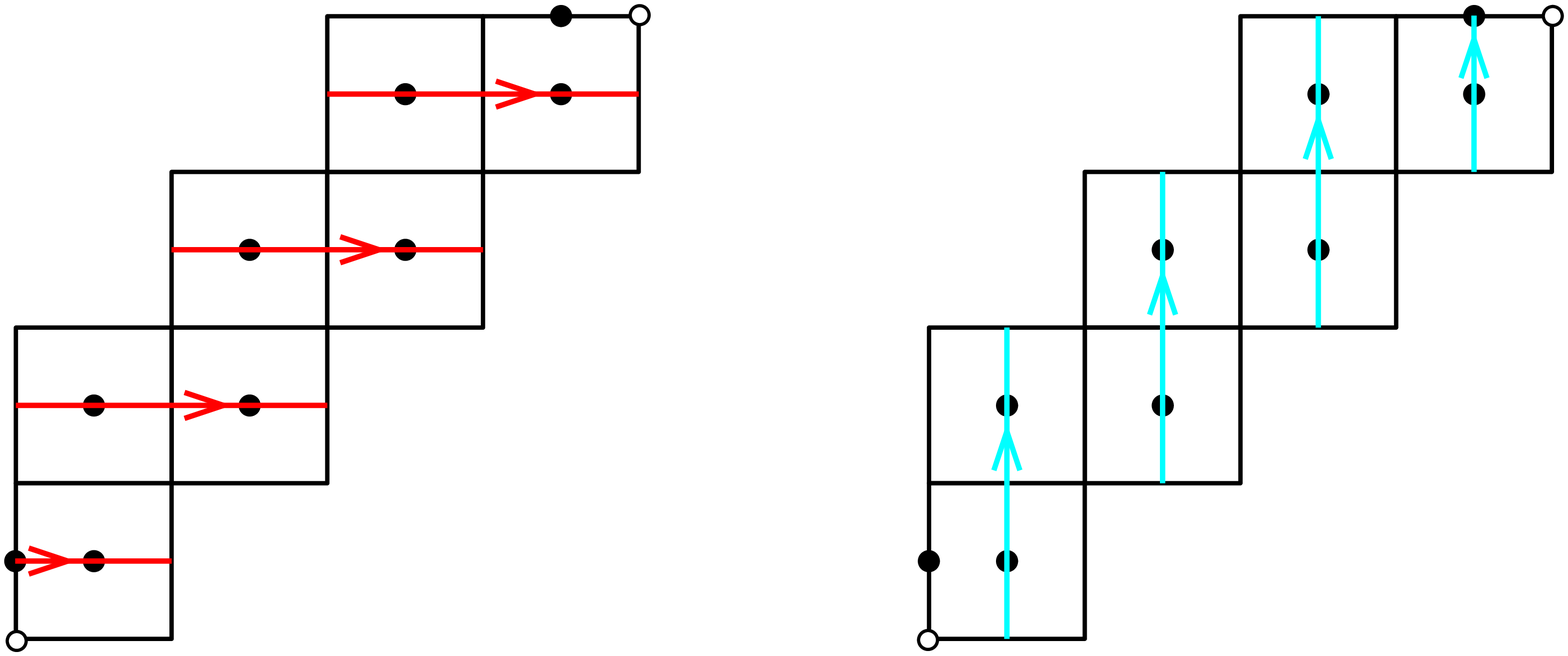} 
\caption{
Simple closed curves $\alpha_1, \dots, \alpha_{g},\beta_1, \dots, \beta_{g}$.}
\label{alpha}
\end{figure}

Let $g\geq 2$. 
In this section, we construct a symplectic basis 
of our algebraic curve $X$ defined by $w^2=z(z^2-1)(z^2-a_1^2)\cdots(z^2-a_{g-1}^2)$
($1<a_1<a_2<\cdots<a_{g-1}$).
By Theorem \ref{representation}, 
the algebraic curve $X$ 
is obtained from the polygon 
\begin{align*}
P=P_0\cup \bigcup_{i=1}^{g-1} \left(P_i\cup Q_i\right) 
\end{align*} 
as in 
Figure \ref{construction_even} or \ref{paper_construction_odd} 
by adjusting the lengths of rectangles. 
Since we construct a symplectic basis of $H_1(X,\mathbb{Z})$ 
by  topological arguments, 
we may describe $P_0$, $P_1$,$\dots, P_{g-1}$, $Q_1$, $\dots, Q_{g-1}$ 
as unit squares.
Let $\alpha_1, \dots, \alpha_{g}, \beta_1, \dots, \beta_{g}$ 
be simple closed curves of $X$ as in Figure \ref{alpha}.
That is, 
$\alpha_1, \dots, \alpha_g$ 
are curves induced by 
horizontal segments passing through centers of the squares 
$P_0$, $P_1$,$\dots, P_{g-1}$, $Q_1$, $\dots, Q_{g-1}$
which are not homotopic to each other. 
We label $\alpha_1, \dots, \alpha_g$ 
so that $\alpha_i$ is below $\alpha_j$
if  $i<j$.
The curves $\beta_1, \dots, \beta_g$ 
are induced by 
vertical segments passing through centers of the squares 
$P_0$, $P_1$,$\dots, P_{g-1}$, $Q_1$, $\dots, Q_{g-1}$
which are not homotopic to each other. 
We label $\beta_1, \dots, \beta_g$ 
so that $\beta_i$ is on the left side of $\beta_j$
if  $i<j$.
We show the following.

\prop{\label{symplectic}
Let  $\alpha_1, \dots, \alpha_{g},\beta_1, \dots, \beta_{g}$ 
be simple closed curves defined as above. 
There exists a symplectic basis 
$\left\{ \alpha_1, \dots, \alpha_{g},  \gamma_1, \dots, \gamma_{g} \right\}$ 
of $H_1(X,\mathbb{Z})$ 
such  that 
$\displaystyle \gamma_j=\sum_{k=j}^g (-1)^{k-j}\beta_k$ 
holds in $H_1(X, \mathbb{Z})$
for all $j=1,2,\dots, g$.
}

\pf{
If $g=2$, we define the simple closed curves $\gamma_1$ and $\gamma_2$ 
as in Figure \ref{g=2}. 
Then, $\gamma_1=\beta_1-\beta_2$ and $\gamma_2=\beta_2$ 
hold in $H_1(X, \mathbb{Z})$ 
and $\left\{\alpha_1, \alpha_2, \gamma_1, \gamma_2 \right\}$ 
is a symplectic basis of  $H_1(X, \mathbb{Z})$.
\begin{figure}[h]
\labellist
\hair 0pt
\pinlabel $\gamma_1$  at  36 202
\pinlabel $\gamma_1$ at 157 202
\pinlabel $\gamma_2$ at 118 202
\endlabellist
\centering
\includegraphics[scale=0.42]{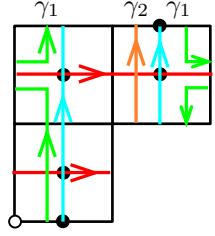} 
\caption{Simple closed curves $\gamma_1$ and $\gamma_2$ if $g=2$.}
\label{g=2}
\end{figure}

Next, we construct simple closed curves 
$\gamma_1$, $\gamma_2$ and $\gamma_3$ 
for the case where $g=3$. 
As in Figure \ref{g=3}, 
we glue a rectangle constructed by two unit squares 
with 
the polygon of Figure \ref{g=2}. 
We reconstruct $\gamma_1$ and $\gamma_2$ 
by gluing curves with the same colors 
as in Figure \ref{g=3} 
and define $\gamma_3=\beta_3$. 
Then $\gamma_1=\beta_1-\beta_2+\beta_3$, 
$\gamma_2=\beta_2-\beta_3$ 
and $\gamma_3=\beta_3$ 
hold  in $H_1(X, \mathbb{Z})$ 
and $\left\{\alpha_1, \alpha_2, \alpha_3, \gamma_1, \gamma_2, \gamma_3 \right\}$ 
is a symplectic basis of  $H_1(X, \mathbb{Z})$ .
\begin{figure}[h]
\labellist
\hair 0pt
\pinlabel $\gamma_1$  at  433 202
\pinlabel $\gamma_1$  at  562 295
\pinlabel $\gamma_2$ at 515 295
\pinlabel $\gamma_2$ at 652 295
\pinlabel $\gamma_3$ at 613 295
\pinlabel $\Rightarrow$ at 326 141
\endlabellist
\centering
\includegraphics[scale=0.42]{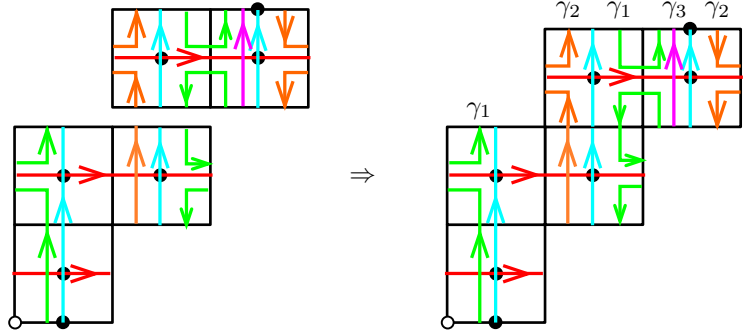} 
\caption{Simple closed curves $\gamma_1$, $\gamma_2$ and $\gamma_3$ if $g=3$}
\label{g=3}
\end{figure}

For the case where $g=4$, 
we glue a rectangle constructed by two unit squares 
with 
the right polygon of Figure \ref{g=3} 
as in Figure \ref{g=4}. 
We reconstruct $\gamma_1$, $\gamma_2$ and $\gamma_3$ 
by gluing curves with the same colors 
as in Figure \ref{g=4} 
and define $\gamma_4=\beta_4$. 
Then $\gamma_1=\beta_1-\beta_2+\beta_3-\beta_4$, 
$\gamma_2=\beta_2-\beta_3+\beta_4$,
 $\gamma_3=\beta_3-\beta_4$ 
 and $\gamma_4=\beta_4$ 
hold  in $H_1(X, \mathbb{Z})$ 
and $\left\{\alpha_1, \alpha_2, \alpha_3, \alpha_4, 
\gamma_1, \gamma_2, \gamma_3, \gamma_4 \right\}$ 
is a symplectic basis of  $H_1(X, \mathbb{Z})$.
\begin{figure}[h]
\labellist
\hair 0pt
\pinlabel $\gamma_1$  at  533 202
\pinlabel $\gamma_1$  at  690 379
\pinlabel $\gamma_2$ at 615 295
\pinlabel $\gamma_2$ at 752 379
\pinlabel $\gamma_3$ at 720 379
\pinlabel $\gamma_3$ at 835 379
\pinlabel $\gamma_4$ at 803 379
\pinlabel $\Rightarrow$ at 411 186
\endlabellist
\centering
\includegraphics[scale=0.42]{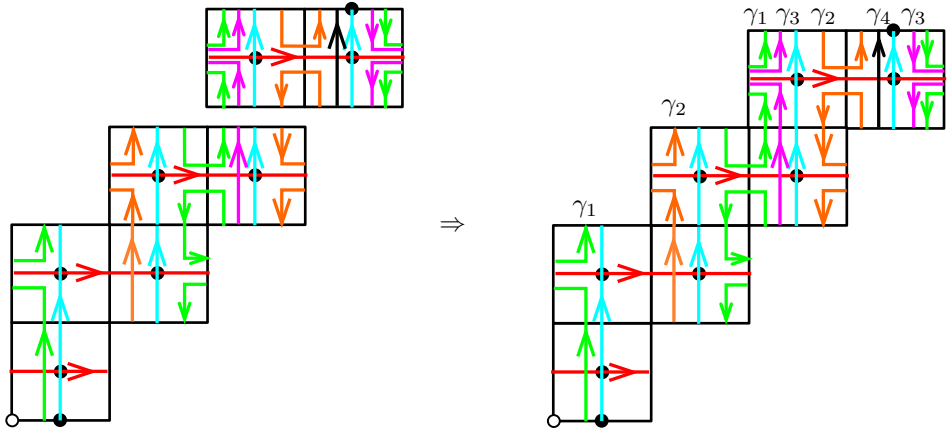} 
\caption{Simple closed curves 
$\gamma_1$, $\gamma_2$, $\gamma_3$ and $\gamma_4$ if $g=4$}
\label{g=4}
\end{figure}

Repeating this process, 
we can construct a symplectic basis 
$\left\{ \alpha_1, \dots, \alpha_{g},  \gamma_1, \dots, \gamma_{g} \right\}$ 
of $H_1(X,\mathbb{Z})$ 
satisfying the condition as in Proposition \ref{symplectic} 
for any $g \geq 2$.
}

\subsection{Period Matrices of $X$}

Let $g\geq 2$ 
and $a_1, a_2, \dots, \alpha_{g-1} \in \mathbb{R}$ 
with $1<a_1<a_2<\cdots<a_{g-1}$. 
Set $f(z)=z(z^2-1)(z^2-a_1^2)\cdots(z^2-a_{g-1}^2)$. 
We calculate the period matrix 
of the algebraic curve 
$X$ defined by $w^2=f(z)$ 
for the symplectic basis 
$\left\{ \alpha_1, \dots, \alpha_{g},  \gamma_1, \dots, \gamma_{g} \right\}$ 
of $H_1(X,\mathbb{Z})$ 
which we construct in Proposition \ref{symplectic}. 

We recall some definitions and theorems.
By Theorem \ref{representation}, 
$X$ is  obtained from the polygon $P$
as in 
Figure \ref{construction_even} 
or Figure \ref{paper_construction_odd}. 
We define some points of $X$ and $P$ 
as Notation\ref{notation}. 
By Theorem \ref{algebraic}, 
we have 
$p_0=(0, 0), p_1=(1, 0), p_{j+1}=(a_j, 0), 
q_1=(-1, 0), q_{j+1}=(-a_j, 0)$ 
for all $j=1, \dots, g-1$. 
Moreover, $X$ has 
an automorphism 
$\tau(z, w)=(-z, i w)$ and 
antiautomorphism $\mu(z, w)=(-\bar{z}, -i\bar{w})$. 
The actions of these maps are seen in subsection \ref{construction_subsection}.
Let $\alpha_1, \dots, \alpha_{g},  \beta_1, \dots, \beta_{g}$ 
be curves defined in subsection \ref{basis} (see Figure \ref{alpha}). 
The curve $\gamma_j$ is defined 
by Proposition \ref{symplectic} 
so that 
$\displaystyle \gamma_j=\sum_{k=j}^g (-1)^{k-j}\beta_k$ 
holds in $H_1(X, \mathbb{Z})$
for all $j=1,2,\dots, g$. 
The family 
$\left\{ \alpha_1, \dots, \alpha_{g},  \gamma_1, \dots, \gamma_{g} \right\}$ 
is a symplectic basis 
of $H_1(X,\mathbb{Z})$.

To calculate the period matrix of $X$ for 
$\left\{ \alpha_1, \dots, \alpha_{g},  \gamma_1, \dots, \gamma_{g} \right\}$, 
we need a basis of  the space $\Omega(X)$ of holomorphic $1$-forms on $X$.
We set $\omega_j=\dfrac{z^{j-1}dz}{w}$ for all $j=1,2,\dots, g$. 
Then $\left\{\omega_1, \omega_2, \dots, \omega_g\right\}$ 
is a basis of $\Omega(X)$. 
Now we calculate the period matrix $\Pi$ 
of $X$ for the symplectic basis 
$\left\{ \alpha_1, \dots, \alpha_{g},  \gamma_1, \dots, \gamma_{g} \right\}$. 
We set 
\begin{align*}
A=\left[\int_{\alpha_k} \frac{z^{j-1}dz}{w} \right],  
B=\left[\int_{\beta_k}  \frac{z^{j-1}dz}{w}  \right],  
C=\left[\int_{\gamma_k}  \frac{z^{j-1}dz}{w}  \right].
\end{align*}

If $g$ is even, 
then we set 
\begin{align*}
I_{j,k}=\left\{
\begin{array}{ll} 
\displaystyle
(-1)^{j-1} \int_{a_{g-2k}}^{a_{g-2k+1}}\dfrac{z^{j-1}dz}{\sqrt{|f(z)|}}
& \left(1\leq k \leq \frac{g}{2}-1\right)\\ 
&\\
\displaystyle
(-1)^{j-1}\int_{1}^{a_1}\dfrac{z^{j-1}dz}{\sqrt{|f(z)|}}
&\left(k=\frac{g}{2}\right)  \\ 
&\\
\displaystyle
\int_{0}^{1}\dfrac{z^{j-1}dz}{\sqrt{|f(z)|}}
&\left(k=\frac{g}{2}+1\right)  \\ 
&\\
\displaystyle
\int_{a_{2k-g-3}}^{a_{2k-g-2}}\dfrac{z^{j-1}dz}{\sqrt{|f(z)|}} 
&
\left(\frac{g}{2}+2 \leq k \leq g\right)
\end{array}
\right.
\end{align*}
for all $j=1,2,\dots, g$.
If $g$ is odd, 
then we set
\begin{align*}
I_{j,k}=\left\{
\begin{array}{ll} 
\displaystyle
(-1)^{j-1}\int_{a_{g-2k}}^{a_{g-2k+1}}\dfrac{z^{j-1}dz}{\sqrt{|f(z)|}}
& \left(1\leq k \leq \frac{g-1}{2}\right)\\ 
&\\
\displaystyle
(-1)^{j-1}\int_{0}^{1}\dfrac{z^{j-1}dz}{\sqrt{|f(z)|}}
&\left(k=\frac{g+1}{2}\right)  \\ 
&\\
\displaystyle
\int_{1}^{a_1}\dfrac{z^{j-1}dz}{\sqrt{|f(z)|}}&
\left(k=\frac{g+3}{2}\right)  \\ 
&\\
\displaystyle
\int_{a_{2k-g-3}}^{a_{2k-g-2}}\dfrac{z^{j-1}dz}{\sqrt{|f(z)|}} 
&
\left(\frac{g+5}{2} \leq k \leq g \right)
\end{array}
\right.
\end{align*}
for all $j=1,2,\dots, g$.
In both cases, we set $\Pi_0=[I_{j,k}]$.
Then, we have the following. 

\lem{\label{a-peri}
The equation 
\begin{align*}
\int_{\alpha_k} \omega_j
=2I_{j,k}
\end{align*}
holds for any $j,k\in\left\{1,2,\dots, g\right\}$.
Moreover, $A=2\Pi_0$ holds.
}

\pf{
We prove the case where $g$ is even. 
We can prove the case where $g$ is odd by the same argument.
Suppose that $g$ is even. 
Let $\varphi:X\to \hat{\mathbb{C}} : (z,w)\to z$ 
be the natural projection.  
Then, we have 
\begin{align*}
\varphi(\alpha_k)=\left\{
\begin{array}{ll} 
\left[-a_{g-2k+1} , -a_{g-2k}\right] &\left(1\leq k \leq \frac{g}{2}-1\right)\\ 
&\\
\left[-a_1 , -1\right] &\left(k=\frac{g}{2}\right)\\
&\\
\left[0,1\right] &\left(k=\frac{g}{2}+1\right) \\
& \\
\left[a_{2k-g-3}, a_{2k-g-2} \right]
&
\left(\frac{g}{2}+2 \leq k \leq g\right)
\end{array}
\right.
\end{align*}
for all $k=1,2, \dots, g$.
If $1\leq k \leq \frac{g}{2}-1$, then 
\begin{align*}
\int_{\alpha_k} \omega_j
=2\int_{\varphi(\alpha_k) }\dfrac{z^{j-1}dz}{\sqrt{|f(z)|}}
=2\int_{-a_{g-2k+1}}^{-a_{g-2k}}\dfrac{z^{j-1}dz}{\sqrt{|f(z)|}}
=2(-1)^{j-1}\int_{a_{g-2k+1}}^{a_{g-2k}}\dfrac{z^{j-1}dz}{\sqrt{|f(z)|}}
=2I_{j,k}
\end{align*}
holds. 
If $k=\frac{g}{2}$, then we have
\begin{align*}
\int_{\alpha_k} \omega_j
=2\int_{\varphi(\alpha_k) }\dfrac{z^{j-1}dz}{\sqrt{|f(z)|}}
=2\int_{-a_1}^{-1}\dfrac{z^{j-1}dz}{\sqrt{|f(z)|}}
=2(-1)^{j-1}\int_{1}^{a_1}\dfrac{z^{j-1}dz}{\sqrt{|f(z)|}}
=2I_{j,k}  .
\end{align*}
Next, if $k=\frac{g}{2}+1$, then 
\begin{align*}
 \int_{\alpha_k}\omega_j
=2\int_{\varphi(\alpha_k) }\dfrac{z^{j-1}dz}{\sqrt{|f(z)|}}
=2\int_{0}^{1}\dfrac{z^{j-1}dz}{\sqrt{|f(z)|}}
=2I_{j,k}
\end{align*}
holds.
Finally, if $\frac{g}{2}+2 \leq k \leq g$, then we have 
\begin{align*}
 \int_{\alpha_k} \omega_j
=2\int_{\varphi(\alpha_k) }\dfrac{z^{j-1}dz}{\sqrt{|f(z)|}}
=2\int_{a_{2k-g-3}}^{a_{2k-g-2}}\dfrac{z^{j-1}dz}{\sqrt{|f(z)|}} 
=2I_{j,k}.
\end{align*}
}

\rem{Let $\varphi:X\to \hat{\mathbb{C}} : (z,w)\to z$ be
the natural projection.  
If $g$ is odd, then we have
\begin{align*}
\varphi(\alpha_k)=\left\{
\begin{array}{ll} 
\left[-a_{g-2k+1} , -a_{g-2k}\right] &\left(1\leq k \leq \frac{g-1}{2}\right)\\ 
&\\
\left[-1 ,0\right] &\left(k=\frac{g+1}{2}\right) \\
&\\
\left[1,a_1\right] &\left(k=\frac{g+3}{2}\right)\\
& \\
\left[a_{2k-g-3}, a_{2k-g-2} \right]
&
\left(\frac{g+5}{2} \leq k \leq g \right)
\end{array}
\right.
\end{align*} 
for all $k=1,2, \dots, g$.
}

Next, we describe the matrix $B$ by $\Pi_0$.

\lem{\label{b-peri}
For any $j,k\in\left\{1,2,\dots, g\right\}$, 
the equation 
\begin{align*}
\int_{\beta_k} \omega_j
=2i(-1)^{j-1}I_{j, g+1-k}
\end{align*}
holds. 
Moreover, we have 
$B
=2i \left[(-1)^{j-1}\delta_{j,k}\right]\Pi_0 \left[\delta_{j+k, g+1}\right]$.
}

\pf{
The automorphism $\tau$ 
maps 
the curve $\alpha_{g+1-k}$ to $\beta_k$ 
for  any $k\in \left\{1,2,\dots, g\right\}$.
Since $\tau(z, w)=(-z, i w)$ holds for all $(z,w)\in X$,  we have 
\begin{align*}
 \tau^\ast \omega_j
 =\tau^\ast \frac{z^{j-1}dz}{w}
=\frac{(-z)^{j-1}(-dz) }{i w}
=i (-1)^{j-1}\frac{z^{j-1}dz}{w}
=i (-1)^{j-1}\omega_j.
\end{align*} 
Therefore, 
by Lemma \ref{a-peri}, 
\begin{align*}
\int_{\beta_k} \omega_j
=\int_{\tau(\alpha_{g+1-k})} \omega_j
=\int_{\alpha_{g+1-k}} \tau^\ast \omega_j
=i (-1)^{j-1} \int_{\alpha_{g+1-k}} \omega_j
=2i(-1)^{j-1}I_{j, g+1-k}
\end{align*}
holds. 
By computation, we obtain 
$B
=i \left[(-1)^{j-1}\delta_{j,k}\right]A \left[\delta_{j+k, g+1}\right]$.
}

Finally, we describe $C$ by $\Pi_0$.

\lem{\label{c-peri}
For any $j,k\in\left\{1,2,\dots, g\right\}$, 
$C
=2i M
\Pi_0 
N
$ holds. 
Here, $M=\Bigl[(-1)^{j-1}\delta_{j,k}\Bigr]$  
and 
$$N=
\left[
\begin{array}{cccccc} 
(-1)^{g-1}  &(-1)^{g}   &(-1)^{g-1}   &\cdots &         -1   &1\\
(-1)^{g} &(-1)^{g-1} & \iddots &    \iddots &\iddots&0\\
 (-1)^{g-1}&\iddots &\iddots  &\iddots&  \iddots    &\vdots\\ 
\vdots& \iddots&  \iddots     &\iddots    &           &\vdots\\ 
-1 & 1 &\iddots  &    &           &\vdots\\ 
1&   0  & \cdots & \cdots  &\cdots &0
\end{array}
\right].
$$
}

\pf{
By the construction of $\gamma_1, \gamma_2, \dots, \gamma_g$ 
and Lemma \ref{b-peri}, 
the equation  
\begin{align*}
 C&=B 
\left[
\begin{array}{cccccc} 
1&   0  & \cdots & \cdots  &\cdots &0\\ 
-1 & 1 &\ddots  &    &           &\vdots\\ 
1& -1&  \ddots     &\ddots    &           &\vdots\\ 
 -1&\ddots &\ddots  &\ddots&  \ddots    &\vdots\\ 
\vdots &\ddots & \ddots &    \ddots &\ddots&0\\
(-1)^{g-1}  &\cdots     &-1  &1 &         -1   &1
\end{array}
\right]\\
&=
i M A 
\left[
\begin{array}{cccccc} 
(-1)^{g-1}  &(-1)^{g}   &(-1)^{g-1}   &\cdots &         -1   &1\\
(-1)^{g} &(-1)^{g-1} & \iddots &    \iddots &\iddots&0\\
 (-1)^{g-1}&\iddots &\iddots  &\iddots&  \iddots    &\vdots\\ 
\vdots& \iddots&  \iddots     &\iddots    &           &\vdots\\ 
-1 & 1 &\iddots  &    &           &\vdots\\ 
1&   0  & \cdots & \cdots  &\cdots &0
\end{array}
\right]
\end{align*}
holds.
}

Now, we have the following theorem.

\thm{\label{period_thm}
Let $g\geq 2$ 
and $a_1, a_2, \dots, \alpha_{g-1} \in \mathbb{R}$ 
with $1<a_1<a_2<\cdots<a_{g-1}$.
Let $\Pi$ be the period matrix of 
$X$ defined by $w^2=z(z^2-1)(z^2-a_1^2)\cdots(z^2-a_{g-1}^2)$ 
for the symplectic basis  
$\left\{ \alpha_1, \dots, \alpha_{g},  \gamma_1, \dots, \gamma_{g} \right\}$ 
constructed in Proposition \ref{symplectic} . 
Then, 
$$
\Pi=i\Pi_0^{-1}M \Pi_0 N 
$$
holds. 
Here, 
$\Pi_0= [I_{j,k}]$, 
$M=\Bigl[(-1)^{j-1}\delta_{j,k}\Bigr]$ 
and 
\begin{align*}
N=
\left[
\begin{array}{cccccc} 
(-1)^{g-1}  &(-1)^{g}   &(-1)^{g-1}   &\cdots &         -1   &1\\
(-1)^{g} &(-1)^{g-1} &  &    \iddots &\iddots&0\\
 (-1)^{g-1}& &\iddots  &\iddots&  \iddots    &\vdots\\ 
\vdots& \iddots&  \iddots     &\iddots    &           &\vdots\\ 
-1 & 1 &\iddots  &    &           &\vdots\\ 
1&   0  & \cdots & \cdots  &\cdots &0
\end{array}
\right].
\end{align*}
Especially,  we have 
$\mathrm{Re}(\Pi)=O$ and  
$\det \left(\Pi \right)=i^g$.
}

\pf{By Lemma \ref{a-peri} and Lemma \ref{c-peri}, 
we have 
\begin{align*}
 \Pi
 =A^{-1}C
 =\frac{1}{2}\Pi_0^{-1} \cdot 
 2i \left[(-1)^{j-1}\delta_{j,k}\right]
\Pi_0 
N
=i\Pi_0^{-1}M \Pi_0 N.
\end{align*}
Since $M$, $N$ and $\Pi_0$ are real matrices, 
$\mathrm{Re}(\Pi)=O$  holds. 
By computation, we obtain 
$\det \left(\Pi \right)=i^g$.
}

\section{Examples }\label{section_example}

In this section, 
we see some examples of period matrices $\Pi$  
of the algebraic curves 
$X$ 
defined by 
$w^2=z(z^2-1)(z^2-a_1^2)\cdots(z^2-a_{g-1}^2)
$ 
for some
 $a_1, a_2, \dots, \alpha_{g-1} \in \mathbb{R}$ 
with $1<a_1<a_2<\cdots<a_{g-1}$.
We calculate them by applying Theorem \ref{period_thm}. 
The definitions of $\Pi_0$, $M$, $N$ are same as in section \ref{section_calculate}.

\subsection{Genus two case}.
If $g=2$, 
then 
the algebraic equation of $X$ 
is of the form 
\begin{align*}
w^2=z(z^2-1)(z^2-a^2)
\end{align*}
where $1<a$.
Set $f(z)=z(z^2-1)(z^2-a^2)$ and 
\begin{align*}
\begin{array}{cc} 
\displaystyle p= \int_{1}^a  \dfrac{dz}{\sqrt{|f(z)|}},& 
\displaystyle  q= \int_{0}^1  \dfrac{dz}{\sqrt{|f(z)|}}, 
\\ 
\displaystyle
r=-\int_{1}^a  \dfrac{z dz}{\sqrt{|f(z)|}},& 
\displaystyle
s= \int_{0}^1  \dfrac{z dz}{\sqrt{|f(z)|}}.
\end{array}
\end{align*}
In this case, we have 
\begin{align*}
\Pi_0=\left[
\begin{array}{ccc} 
p&q \\ 
r &s 
\end{array}
\right], 
M=\left[
\begin{array}{ccc} 
1&0 \\ 
0 &-1 
\end{array}
\right],
N=\left[
\begin{array}{cc} 
-1&1 \\  
1 &0 
\end{array}
\right]. 
\end{align*}
By Theorem \ref{period_thm}, 
the period matrix $\Pi$ of $X$ is 
described as 
$$
\Pi 
=\Pi_0^{-1}M\Pi_0 N 
=
\frac{i}{ps-qr}
\left[
\begin{array}{cc} 
2qs-ps-qr& ps+qr \\ 
2pr-ps-qr&  -2pr\\ 
\end{array}
\right].
$$
Since the period matrix 
$\Pi$
is symmetric, 
we have 
$\Pi_{1, 2}=\Pi_{2,1}=\frac{1}{2}\left(\Pi_{1, 2}+\Pi_{2,1}\right)=pr$.
and 
$pr-ps-qr=0$.
Therefore, we have the following.
\thm{
Let  $a>1$ and $X$ 
 an algebraic curve defined 
by  
$
w^2=z(z^2-1)(z^2-a^2)
$.
Set $f(z)=z(z^2-1)(z^2-a^2)$ 
and 
$p$, $q$, $r$, $s$ as above. 
Then, the period matrix of $X$ is 
$$
\Pi 
=
\frac{i}{ps-qr}
\left[
\begin{array}{cc} 
2qs-pr& pr \\ 
pr&  -2pr\\ 
\end{array}
\right].
$$
}

\exam
{The following are examples calculated by Mathematica. 
\begin{enumerate}
\item For the algebraic curve 
$X$ 
defined by 
$w^2=z(z^2-1)(z^2-2)$, 
the period matrix $\Pi$ satisfies 
\begin{align*}
\Pi
 \fallingdotseq
i\left[
\begin{array}{cc}
 1.42594 & -0.409423 \\
 -0.409423 & 0.818846 \\
\end{array}
\right]
\end{align*}

\item For the algebraic curve 
$X$ 
defined by 
$w^2=z(z^2-1)(z^2-4)$, 
the period matrix $\Pi$ satisfies 
\begin{align*}
\Pi
 \fallingdotseq
 i\left[
\begin{array}{cc}
 1.25352 & -0.497668 \\
 -0.497668 & 0.995336 \\
\end{array}
\right].
\end{align*}

\item For the algebraic curve 
$X$ 
defined by 
$w^2=z(z^2-1)(z^2-1.0001^2)$, 
the period matrix $\Pi$ satisfies 
\begin{align*}
\Pi
 \fallingdotseq
i\left[
\begin{array}{cc}
 3.87984 & -0.131086 \\
 -0.131086 & 0.262171 \\
\end{array}
\right].
\end{align*}
\end{enumerate}

}

\subsection{Genus three case}
If $g=3$, 
then 
the algebraic equation of $X$ 
is of the form 
\begin{align*}
w^2=z(z^2-1)(z^2-a^2)(z^2-b^2)
\end{align*}
where $1<a<b$.
Set $f(z)=z(z^2-1)(z^2-a^2)(z^2-b^2)$ and 
\begin{align*}
\begin{array}{ccc} 
\displaystyle I_{1,1}= \int_{a}^b  \dfrac{dz}{\sqrt{|f(z)|}},& 
\displaystyle I_{1,2}= \int_{0}^1  \dfrac{dz}{\sqrt{|f(z)|}},&
 \displaystyle  I_{1,3}= \int_{1}^a  \dfrac{dz}{\sqrt{|f(z)|}} ,\\ 
&\\
\displaystyle I_{2,1}= -\int_{a}^b  \dfrac{zdz}{\sqrt{|f(z)|}},& 
\displaystyle I_{2,2}= -\int_{0}^1  \dfrac{zdz}{\sqrt{|f(z)|}},&
 \displaystyle  I_{2,3}= \int_{1}^a  \dfrac{zdz}{\sqrt{|f(z)|}}, \\ 
&\\
\displaystyle I_{3,1}= \int_{a}^b  \dfrac{z^2dz}{\sqrt{|f(z)|}},& 
\displaystyle I_{3,2}= \int_{0}^1  \dfrac{z^2dz}{\sqrt{|f(z)|}},&
 \displaystyle  I_{3,3}= \int_{1}^a  \dfrac{z^2dz}{\sqrt{|f(z)|}}.
\end{array}
\end{align*}

\begin{align*}
\Pi_0=[I_{j,k}], 
M=\left[
\begin{array}{ccc} 
1&0 &0\\  
0 &-1 &0\\
0 & 0 & 1 
\end{array}
\right],
N=\left[
\begin{array}{ccc} 
1 & -1 & 1\\
-1&1  & 0\\ 
1 &0  & 0
\end{array}
\right].
\end{align*}
By Theorem \ref{period_thm}, 
the period matrix $\Pi$ of $X$ is 
described as 
$
\Pi 
=\Pi_0^{-1}M\Pi_0 N$.

\exam
{The following are examples calculated by Mathematica. 
\begin{enumerate}
\item 
For the algebraic curve 
$X$ 
defined by 
$w^2=x(z^2-1)(z^2-4)(z^2-9)$, 
the period matrix $\Pi$ satisfies 
\begin{align*}
\Pi_0 
 \fallingdotseq
i
\left[
\begin{array}{ccc}
 1.39658 & -0.687212 & 0.371981 \\
 -0.687212 & 1.2467 & -0.495331 \\
 0.371981 & -0.495331 & 0.994534 \\
\end{array}
\right].
\end{align*}

\item 
For the algebraic curve 
$X$ 
defined by 
$w^2=z(z^2-1)(z^2-4)(z^2-10000)$, 
the period matrix $\Pi$ satisfies 
\begin{align*}
 \Pi
 \fallingdotseq
i
\left[
\begin{array}{ccc}
 1.0086 & -0.915883 & 0.869095 \\
 -0.915883 & 1.82283 & -1.28051 \\
 0.869095 & -1.28051 & 1.99277 \\
\end{array}
\right].
\end{align*}

\item 
For the algebraic curve 
$X$ 
defined by 
$w^2=z(z^2-1)(z^2-1.00001^2)(z^2-1.0001^2)$, 
the period matrix $\Pi$ satisfies 
\begin{align*}
 \Pi
 \fallingdotseq
i
\left[
\begin{array}{ccc}
 1.61889 & -0.996731 & 0.0052455 \\
 -0.996731 & 1.00002 & -0.00525414 \\
 0.00524596 & -0.00525459 & 1.59888 \\
\end{array}
\right].
\end{align*}
\end{enumerate}
}

\subsection{genus four case}
If $g=4$, 
then 
the algebraic equation of $X$ 
is of the form 
\begin{align*}
w^2=z(z^2-1)(z^2-a^2)(z^2-b^2)
\end{align*}
where $1<a<b<c$.
Set $f(z)=z(z^2-1)(z^2-a^2)(z^2-b^2)$ and 
\begin{align*}
\begin{array}{cccc} 
\displaystyle I_{1,1}= \int_{b}^c  \dfrac{dz}{\sqrt{|f(z)|}},& 
\displaystyle I_{1,2}= \int_{1}^a  \dfrac{dz}{\sqrt{|f(z)|}},&
 \displaystyle  I_{1,3}= \int_{0}^1  \dfrac{dz}{\sqrt{|f(z)|}}, &
  \displaystyle  I_{1,4}= \int_{a}^b  \dfrac{dz}{\sqrt{|f(z)|}},
 \\ 
&\\
\displaystyle I_{2,1}= -\int_{b}^c  \dfrac{z dz}{\sqrt{|f(z)|}},& 
\displaystyle I_{2,2}= -\int_{1}^a  \dfrac{z dz}{\sqrt{|f(z)|}},&
 \displaystyle  I_{2,3}= \int_{0}^1  \dfrac{z dz}{\sqrt{|f(z)|}},&
  \displaystyle  I_{2,4}= \int_{a}^b  \dfrac{z dz}{\sqrt{|f(z)|}},
 \\ 
&\\
\displaystyle I_{3,1}= \int_{b}^c  \dfrac{z^2 dz}{\sqrt{|f(z)|}},& 
\displaystyle I_{3,2}= \int_{1}^a  \dfrac{z^2 dz}{\sqrt{|f(z)|}},&
 \displaystyle  I_{3,3}= \int_{0}^1  \dfrac{z^2 dz}{\sqrt{|f(z)|}}, &
  \displaystyle  I_{3,4}= \int_{a}^b  \dfrac{z^2 dz}{\sqrt{|f(z)|}} ,
  \\ 
&\\
\displaystyle I_{4,1}= -\int_{b}^c  \dfrac{z^3 dz}{\sqrt{|f(z)|}},& 
\displaystyle I_{4,2}= -\int_{1}^a  \dfrac{z^3 dz}{\sqrt{|f(z)|}},&
 \displaystyle  I_{4,3}= \int_{0}^1  \dfrac{z^3 dz}{\sqrt{|f(z)|}}, &
  \displaystyle  I_{4,4}= \int_{a}^b  \dfrac{z^3 dz}{\sqrt{|f(z)|}} .
\end{array}
\end{align*}
In this case, we have 
\begin{align*}
\Pi_0=[I_{j,k}], 
M=\left[
\begin{array}{cccc} 
1&0 &0 & 0\\  
0 &-1 &0 & 0\\
0 & 0 & 1  &0\\
0 & 0 & 0 &-1
\end{array}
\right],
N=\left[
\begin{array}{cccc} 
-1& 1 & -1 &1\\
1 & -1 & 1 & 0\\
-1&1  & 0 &0\\ 
1 &0  & 0 &0
\end{array}
\right].
\end{align*}
By Theorem \ref{period_thm}, 
the period matrix $\Pi$ of $X$ is 
described as 
$
\Pi 
=\Pi_0^{-1}M\Pi_0 N$.

\exam
{The following are examples calculated by Mathematica. 
\begin{enumerate}
\item 
For the algebraic curve 
$X$ 
defined by 
$ w^2=z(z^2-1)(z^2-4)(z^2-9)(z^2-16)$, 
the period matrix $\Pi$ satisfies 
\begin{align*}
\Pi
 \fallingdotseq
i
\left[
\begin{array}{cccc}
 1.49592 & -0.805976 & 0.529694 & -0.309252 \\
 -0.805976 & 1.3887 & -0.683972 & 0.370541 \\
 0.529694 & -0.683972 & 1.24537 & -0.494738 \\
 -0.309252 & 0.370541 & -0.494738 & 0.99427 \\
\end{array}
\right].
\end{align*}

\item 
For the algebraic curve 
$X$ 
defined by 
$w^2=z(z^2-1)(z^2-1.00001^2)(z^2-1.0001^2)(z^2-1.001^2)$, 
the period matrix $\Pi$ satisfies 
\begin{align*}
\Pi 
 \fallingdotseq
i
\left[
\begin{array}{cccc}
 3.19594 & -2.99265 & 0.161826 & -0.141682 \\
 -2.99265 & 4.42015 & -0.161842 & 0.141696 \\
 0.161828 & -0.161844 & 0.323658 & -0.283311 \\
 -0.141685 & 0.141699 & -0.283311 & 0.861195 \\
\end{array}
\right]. 
\end{align*}

\item 
For the algebraic curve 
$X$ 
defined by 
$w^2=z(z^2-1)(z^2-1.001^2)(z^2-1.01^2)(z^2-100^2)$, 
the period matrix $\Pi$ satisfies 
\begin{align*}
\Pi 
 \fallingdotseq
i
\left[
\begin{array}{cccc}
 1.00423 & -0.996561 & 0.912285 & -0.9105 \\
 -0.996562 & 2.58767 & -0.957418 & 0.954779 \\
 0.912286 & -0.957419 & 1.82628 & -1.79174 \\
 -0.9105 & 0.95478 & -1.79174 & 2.38377 \\
\end{array}
\right].
\end{align*}

\item 
For the algebraic curve 
$X$ 
defined by 
$w^2=z(z^2-1)(z^2-1.001^2)(z^2-1000^2)(z^2-10000000^2)$, 
the period matrix $\Pi$ satisfies 
\begin{align*}
\Pi 
 \fallingdotseq
i
\left[
\begin{array}{cccc}
 1.00004 & -0.99103 & 0.990938 & -0.990902 \\
 -0.991031 & 5.08919 & -2.10865 & 1.95408 \\
 0.990937 & -2.10865 & 2.29095 & -1.9742 \\
 -0.990902 & 1.95408 & -1.9742 & 1.98199 \\
\end{array}
\right].
\end{align*}

\end{enumerate}
}

%
%
%
%

%
%

\section{Appendix}\label{section_appendix}

In this section, 
we show 
that our algebraic curves 
are not conformal equivalent to 
four 
algebraic curves 
whose period matrices are calculated by 
Schindler\cite{Schindler93}. 
That is, we show the following.

\thm{ 
Let $g\geq 2$ and
$a_1, a_2, \dots, a_{g-1} \in \mathbb{R}$ 
with $1<a_1<a_2<\cdots<a_{g-1}$. 
Let $X$ be the algebraic curve defined by 
$\displaystyle w^2=z(z^2-1)(z^2-a_1^2)\cdots(z^2-a_{g-1}^2)$.
Then, $X$ 
is not conformal equivalent to 
algebraic curves 
$Y_1$, $Y_2$, $Y_3$, $Y_4$ 
defined by the algebraic equations  
$w^2=z^{2g+2}-1$, 
$w^2=z(z^{2g+1}-1)$, 
$w^2=z(z^{2g}-1)$ and 
$w^2=z^{2g+1}-1$, respectively. 
}

%
 

This theorem is obtained from 
Lemma \ref{alg1}, 
\ref{alg2},  
\ref{alg3} 
and 
\ref{alg4}.

\lem{\label{alg1}
The algebraic curve $X$ is not conformal 
equivalent to 
the algebraic curve 
$Y_1$ defined by $w^2=z^{2g+2}-1$.
}

\pf{
Assume that there exists 
a conformal map $\widetilde{\varphi} : X \to Y_1$.
Denote by $\sigma_X$ and $\sigma_{Y_1}$ 
the hyperelliptic involutions of $X$ and $Y_1$, respectively.
Let $p_X: X\to \hat{\mathbb{C}} ; (z,w)\mapsto z$ 
and $p_{Y_1}: Y_1\to \hat{\mathbb{C}} ; (z,w)\mapsto z$ 
be natural projections. 
Since 
$\sigma_{Y_1}\circ \widetilde{\varphi}=\widetilde{\varphi} \circ \sigma_X$
holds, the map 
 $\widetilde{\varphi}$ induces a conformal map 
$\varphi :  \hat{\mathbb{C}} \to  \hat{\mathbb{C}} $ 
satisfying  $p_2\circ \widetilde{\varphi} =\varphi \circ p_1$. 
The map $\varphi$ is a M\"{o}bius transformation. 
Since $\varphi$ maps 
the set $\left\{ 0,\pm 1, \pm a_1, \pm a_2, \dots, \pm a_{g-1}, \infty \right\}$ 
to $\left\{ z: z^{2g+2}=1 \right\}$,  
$\varphi(\hat{\mathbb{R}})=S^1$ holds. 
By composing a rotation about $0$, 
we may assume that $\varphi(0)=1$.
Then, 
one of the following holds; 
\begin{itemize}
\item 
$\varphi(\infty)=-1$, 
$\varphi(1)=e^{\frac{\pi}{g+1}i}$, 
$\varphi(-1)=e^{-\frac{\pi}{g+1}i}$, 
$\varphi(\infty)=-1$, 

\item 
$\varphi(\infty)=-1$, 
$\varphi(1)=e^{-\frac{\pi}{g+1}i}$, 
$\varphi(-1)=e^{\frac{\pi}{g+1}i}$,
$\varphi(\infty)=-1$.
\end{itemize}
By considering the cross ratios, 
we conclude that there exist 
no such M\"{o}bius transformations. 
}

\lem{\label{alg2}
The algebraic curve $X$ is not conformal 
equivalent to 
the algebraic curves 
$Y_2$ defined by $w^2=z(z^{2g+1}-1)$.
}

\pf{
Assume that there exists 
a conformal map $\widetilde{\varphi} : X \to Y_2$.
Denote by $\sigma_X$ and $\sigma_{Y_2}$ 
the hyperelliptic involutions of $X$ and $Y_2$, respectively.
Let $p_X: X\to \hat{\mathbb{C}} ; (z,w)\mapsto z$ 
and $p_{Y_2}: Y_2\to \hat{\mathbb{C}} ; (z,w)\mapsto z$ 
be natural projections. 
Since 
$\sigma_{Y_2}\circ \widetilde{\varphi}=\widetilde{\varphi} \circ \sigma_X$
holds, the map 
 $\widetilde{\varphi}$ induces a conformal map 
$\varphi :  \hat{\mathbb{C}} \to  \hat{\mathbb{C}} $. 
The map $\varphi$ is a M\"{o}bius transformation 
and $\varphi$  maps 
the set $\left\{ 0,\pm 1, \pm a_1, \pm a_2, \dots, \pm a_{g-1}, \infty \right\}$ 
to $\left\{ z: z^{2g+1}=1 \right\} \cup \left\{0\right\}$.   
This contradicts that $\varphi(\hat{\mathbb{R}})$ is a line or a circle. 
}

\lem{\label{alg3}
The algebraic curve $X$ is not conformal 
equivalent to 
the algebraic curves 
$Y_3$ defined by $w^2=z(z^{2g}-1)$.
}

\pf{
Assume that there exists 
a conformal map $\widetilde{\varphi} : X \to Y_3$.
Denote by $\sigma_X$ and $\sigma_{Y_3}$ 
the hyperelliptic involutions of $X$ and $Y_3$, respectively.
Let $p_X: X\to \hat{\mathbb{C}} ; (z,w)\mapsto z$ 
and $p_{Y_3}: Y_3\to \hat{\mathbb{C}} ; (z,w)\mapsto z$ 
be natural projections. 
Since 
$\sigma_{Y_3}\circ \widetilde{\varphi}=\widetilde{\varphi} \circ \sigma_X$
holds, the map 
 $\widetilde{\varphi}$ induces a conformal map 
$\varphi :  \hat{\mathbb{C}} \to  \hat{\mathbb{C}} $. 
The map $\varphi$ is a M\"{o}bius transformation 
and $\varphi$  maps 
the set $\left\{ 0,\pm 1, \pm a_1, \pm a_2, \dots, \pm a_{g-1}, \infty \right\}$ 
to $\left\{ z: z^{2g}=1 \right\} \cup \left\{0, \infty \right\}$.   
This contradicts that $\varphi(\hat{\mathbb{R}})$ is a line or a circle. 
}

\lem{\label{alg4}
The algebraic curve $X$ is not conformal 
equivalent to 
the algebraic curves 
$Y_4$ defined by $w^2=z^{2g+1}-1$.
}

\pf{
Assume that there exists 
a conformal map $\widetilde{\varphi} : X \to Y_4$.
Denote by $\sigma_X$ and $\sigma_{Y_4}$ 
the hyperelliptic involutions of $X$ and $Y_4$, respectively.
Let $p_X: X\to \hat{\mathbb{C}} ; (z,w)\mapsto z$ 
and $p_{Y_4}: Y_4\to \hat{\mathbb{C}} ; (z,w)\mapsto z$ 
be natural projections. 
Since 
$\sigma_{Y_4}\circ \widetilde{\varphi}=\widetilde{\varphi} \circ \sigma_X$
holds, the map 
 $\widetilde{\varphi}$ induces a conformal map 
$\varphi :  \hat{\mathbb{C}} \to  \hat{\mathbb{C}} $. 
The map $\varphi$ is a M\"{o}bius transformation 
and $\varphi$  maps 
the set $\left\{ 0,\pm 1, \pm a_1, \pm a_2, \dots, \pm a_{g-1}, \infty \right\}$ 
to $\left\{ z: z^{2g+1}=1 \right\} \cup \left\{\infty \right\}$.   
This contradicts that $\varphi(\hat{\mathbb{R}})$ is a line or a circle. 
}

%
%
%
\bibliographystyle{alpha}
%

%

\bibliography{ref}

\end{document}